# SYMÉTRIES BIRATIONNELLES DES SURFACES FEUILLETÉES.

SERGE CANTAT ET CHARLES FAVRE

ABSTRACT. We provide a classification of complex projective surfaces with a holomorphic foliation whose group of birational symetries is infinite.

## 1. INTRODUCTION

Depuis quelques années, l'étude des feuilletages holomorphes singuliers des surfaces complexes compactes a connu un essor remarquable, en particulier pour les surfaces projectives [5], [26], [27]. Les travaux de M. Brunella, M. McQuillan et L. G. Mendès permettent ainsi de réaliser une classification de type Enriques-Kodaira pour les surfaces feuilletées [7]. Parallèlement, la dynamique des transformations holomorphes ou méromorphes des variétés kählériennes, notamment des surfaces, s'est rapidement développée [23], [3], [14], [29], [22], [4], [12], [10]. Cet article se place à la croisée de ces deux thèmes. Son but, en effet, est de classer les triplets $(X, \varphi, \mathcal{F})$ formés d'une surface complexe compacte projective $X$, d'une application birationnelle $\varphi : X \circlearrowleft$ d'ordre infini et d'un feuilletage holomorphe $\varphi$-invariant $\mathcal{F}$, éventuellement singulier. Une application birationnelle préserve un feuilletage lorsqu'elle envoie feuille sur feuille là où elle est définie. Nous noterons $\text{Bir}(\mathcal{F})$ le groupe des transformations birationnelles laissant $\mathcal{F}$ invariant et $\text{Aut}(\mathcal{F})$ le sous-groupe formé des transformations biholomorphes.

Pour la majeure partie des feuilletages, $\text{Bir}(\mathcal{F})$ coïncide avec $\text{Aut}(\mathcal{F})$ et est un groupe fini ; le résultat principal de ce texte fournit la liste des feuilletages qui violent ce "principe général" (théorèmes 1.1 et 1.2). De manière analogue, la plupart des applications biméromorphes préservant un feuilletage préservent en réalité un pinceau de courbes rationnelles ou elliptiques. Dans ce cas, leur dynamique se ramène au cas unidimensionnel. Nous donnerons la liste des tranformations birationnelles dont la dynamique n'est pas triviale et qui préservent un feuilletage holomorphe avec une feuille transcendante (voir corollaire 1.3). Nous verrons par exemple que l'existence d'un tel feuilletage assure l'existence d'un second feuilletage invariant : la dynamique de la transformation est alors de type "pseudo-Anosov".

1.1. **Les feuilletages et leurs symétries.** Soient $X$ une surface complexe compacte et $\mathcal{F}$ un feuilletage holomorphe sur $X$, donné par une famille $v_i$ de champs de vecteurs holomorphes à zéros isolés. Les champs $v_i$ sont définis sur les ouverts $\mathcal{U}_i$ d'un recouvrement de $X$ et sont soumis à des relations de compatibilités

$$v_i = g_{ij} v_j, \quad g_{ij} \in \mathcal{O}^*(\mathcal{U}_i \cap \mathcal{U}_j),$$

assurant que les courbes intégrales de $v_i$ et $v_j$ se recollent sur $\mathcal{U}_i \cap \mathcal{U}_j$ pour former les feuilles de $\mathcal{F}$. Le fibré en droites associé au cocycle $\{g_{ij}\}$ ne dépend que de $\mathcal{F}$ : c'est *le fibré cotangent* de $\mathcal{F}$. Il est noté $T^*_\mathcal{F}$ et sera parfois dénommé *fibré canonique*. À ce fibré est associée une dimension de Kodaira $\mathrm{kod}(\mathcal{F}) \in \{-\infty, 0, 1, 2\}$ qui mesure son degré de positivité (voir [7]).

Lorsque $\mathrm{kod}(\mathcal{F})$ est égale à 2, le feuilletage $\mathcal{F}$ est dit *de type général*. C'est le cas pour un feuilletage de $\mathbb{P}^2$ générique. Dans [28], J.-V. Pereira et P.-V. Sanchez ont étendu le théorème d'Andreotti [1] au cas feuilleté en démontrant que $\mathrm{Aut}(\mathcal{F})$ est un groupe fini pour tout feuilletage de type général.

M. Brunella a montré que les feuilletages pour lesquels $\mathrm{Bir}(\mathcal{F})$ contient strictement $\mathrm{Aut}(\mathcal{F})$ sont tous, à une exception près, des fibrations rationnelles ou des feuilletages de Riccati, c'est-à-dire qu'ils sont transverses aux fibres génériques d'une fibration rationnelle [6]. En supposant l'inclusion stricte $\mathrm{Aut}(\mathcal{F}) \subsetneq \mathrm{Bir}(\mathcal{F})$ vérifiée pour tout modèle birationnel de $\mathcal{F}$, nous montrons que le feuilletage est en fait une fibration ou un feuilletage linéaire de $\mathbb{P}^1 \times \mathbb{P}^1$ (théorème 1.2). Il est important de noter que, contrairement au cas des surfaces, cette hypothèse n'implique pas que la dimension de Kodaira (numérique) de $\mathcal{F}$ vaut $-\infty$.

Pour chaque transformation birationnelle $\varphi$ d'une surface complexe compacte, nous noterons $\varphi^*$ l'application linéaire induite par $\varphi$ sur le second groupe de cohomologie de $X$ et $\|\varphi^*\|$ la norme de cette application vis-à-vis d'une métrique arbitraire sur $\mathrm{H}^2(X, \mathbb{R})$. La complexité de la dynamique de $\varphi$ se mesure à l'aide du comportement asymptotique de la suite $\|(\varphi^*)^n\|$ : cette suite peut être bornée, croître linéairement, croître quadratiquement ou croître exponentiellement vite [10], [12]. Nous parlerons de *croissance des degrés* pour décrire ce comportement. Les automorphismes avec une croissance des degrés exponentielle sont exactement ceux d'entropie topologique strictement positive. La croissance des degrés d'un élément générique du groupe des applications birationnelles de $\mathbb{P}^2$ est exponentielle ; c'est le cas, par exemple, des transformations de Hénon. Dans [6], M. Brunella a montré que les transformations de Hénon ne préservent pas de feuilletage holomorphe.

1.2. **Exemples.** Avant d'énoncer précisément nos résultats, il convient de présenter les exemples principaux. Le premier d'entre eux est fourni par les feuilletages stables et instables d'un difféomorphisme holomorphe de type Anosov. Un théorème d'É. Ghys montre que ces difféomorphismes sont donnés par des transformations affines des tores [16].

**Exemple 1.1. a.-** Si $T = \mathbb{C}^2/\Lambda$ est un tore complexe de dimension 2 défini par le réseau $\Lambda \subset \mathbb{C}^2$, toute application affine $\varphi$ préservant $\Lambda$ induit un automorphisme de $T$. Lorsque la partie linéaire de $\varphi$ est d'ordre infini, elle peut être hyperbolique, et $\varphi$ induit alors un automorphisme Anosov préservant deux feuilletages linéaires, ou unipotente, et, dans ce cas, $\varphi$ préserve une fibration elliptique. Les applications $\varphi_1(z,w) = (z+2w, z+w)$ et $\varphi_2(z,w) = (z, z+w)$ sur le produit $\mathbb{C}/\Lambda_0 \times \mathbb{C}/\Lambda_0$ d'une courbe elliptique par elle-même en sont les exemples les plus simples.

**b.-** Parfois, le tore $T$ possède un groupe fini d'automorphismes $G$ normalisé par $\varphi$. L'automorphisme $\psi$ induit sur la désingularisée $X$ du quotient $T/G$ préserve alors les projetés des feuilletages stables et instables de $\varphi$ (lorsque $\varphi$ est Anosov), ou une fibration elliptique. Lorsque $G$ est réduit à l'involution $\sigma(x,y) = (-x,-y)$, $X$ est une surface de Kummer. Par analogie, pour $G$ quelconque, nous dirons que $T/G$ est une surface de Kummer généralisée.

**Exemple 1.2.** Quelques applications biméromorphes d'ordre infini préservent des fibrations elliptiques ou rationnelles. C'est le cas des flots de champ de vecteurs parallèles à de telles fibrations ou de l'automorphisme $\varphi_2$ défini ci-dessus. Enfin, si $X$ est le produit de la droite projective $\mathbb{P}^1$ par une surface de Riemann compacte $B$ et si $\alpha : B \to \overline{\mathbb{C}}$ est une application méromorphe non constante, la transformation $\varphi(z,w) = (\alpha(w) \cdot z, w)$ est une transformation birationnelle de $X$ d'ordre infini qui préserve la fibration rationnelle $X \to B$.

**Exemple 1.3.** Cet exemple est analogue à l'exemple 1.1 mais conduit à des transformations biméromorphes.

**a.-** Toute application rationnelle du type $\varphi(z,w) = (z^a w^b, z^c w^d)$ sur $\mathbb{P}^1 \times \mathbb{P}^1$ sera appelée *transformation monomiale* ; lorsque $|ad - bc| = 1$, c'est une application birationnelle. Supposons que le rayon spectral de la matrice associée,

$$M := \begin{bmatrix} a & b \\ c & d \end{bmatrix},$$

soit de module strictement plus grand que 1. L'application $\varphi$ préserve alors les deux feuilletages holomorphes définis par les 1-formes $\alpha w \, dz + \beta z \, dw$ où $(\alpha, \beta)$ est un vecteur propre de ${}^t M$. Le tore réel $\{|z| = |w| = 1\}$ est invariant par cette transformation ; elle y induit un difféomorphisme Anosov.

**b.-** Le quotient de $\mathbb{P}^1 \times \mathbb{P}^1$ par l'involution $\tau(z,w) = (1/z, 1/w)$ est une surface rationnelle possédant quatre points singuliers. Notons $Y$ la surface obtenue après leur résolution. Les applications monomiales de l'exemple précédent commutent à $\tau$ et induisent donc des applications birationnelles de $Y$ préservant les images des deux feuilletages préservés par $\varphi$.

1.3. **Résultats principaux.** Nous démontrerons les deux théorèmes suivants :

**Théorème 1.1.** *Soit $\mathcal{F}$ un feuilletage holomorphe singulier sur une surface projective $X$ possédant un groupe de symétries holomorphes $\mathrm{Aut}(\mathcal{F})$ infini. Il existe alors au moins un élément $\varphi$ d'ordre infini dans $\mathrm{Aut}(\mathcal{F})$ et la situation est l'une des trois suivantes :*

   **1.-** *$\mathcal{F}$ est invariant par un champ de vecteurs holomorphe ;*

   **2.-** *$\mathcal{F}$ est une fibration elliptique ;*

   **3.-** *la surface $X$ est une surface de Kummer généralisée, $\varphi$ se relève en un automorphisme Anosov $\psi$ du tore et $\mathcal{F}$ est la projection sur $X$ du feuilletage stable ou instable de $\psi$ (exemples 1.1).*

Notons que les deux derniers cas sont mutuellement exclusifs. Lorsque $X$ n'est pas un tore, les feuilletages apparaissant au cas (3.-) ne sont jamais préservés par un champ de vecteurs. Nous décrirons les feuilletages invariants par un champ de vecteurs à la proposition 3.8. Notons enfin qu'il existe des fibrations elliptiques dont le groupe d'automorphismes est fini.

**Théorème 1.2.** *Soit $\mathcal{F}$ un feuilletage tel que l'inclusion stricte $\mathrm{Aut}(\mathcal{F}) \subsetneq \mathrm{Bir}(\mathcal{F})$ soit vérifiée pour* tout *modèle birationnel de $\mathcal{F}$. Alors $\mathrm{Bir}(\mathcal{F})$ possède un élément birationnel d'ordre infini $\varphi$ et*

   **1.-** *soit $\mathcal{F}$ est une fibration rationnelle ;*

   **2.-** *soit la situation est birationnellement conjuguée à celle de l'exemple 1.3.*

Réciproquement, toute fibration rationnelle possède une symétrie birationnelle d'ordre infini qui n'est pas birationnellement conjuguée à un automorphisme.

Nous renvoyons à la section 7 pour quelques corollaires de cette classification ; citons dès à présent le

**Corollaire 1.3.** *Soit $\varphi : X \circlearrowleft$ une application birationnelle d'une surface projective préservant un feuilletage holomorphe $\mathcal{F}$. Si la croissance des degrés de $\varphi$ est exponentielle, il existe un changement de variables birationnel, puis un revêtement fini par un tore (resp. par $\mathbb{P}^1 \times \mathbb{P}^1$), tel que $\varphi$ se relève en un automorphisme Anosov du tore (resp. une application monomiale de $\mathbb{P}^1 \times \mathbb{P}^1$) et le feuilletage $\mathcal{F}$ en un feuilletage linéaire. En particulier, $\varphi$ préserve automatiquement deux feuilletages.*

1.4. **Organisation du texte.** La preuve du théorème 1.1 se fait en plusieurs étapes. Nous montrerons à la proposition 3.9 que, si le groupe $\mathrm{Aut}(\mathcal{F})$ est infini, il possède un élément $\varphi$ d'ordre infini. Il s'agit alors d'étudier les feuilletages invariants par un automorphisme $\varphi$ d'ordre infini.

Lorsque la croissance des degrés de $\varphi$ est exponentielle, nous montrons que la dimension de Kodaira de $\mathcal{F}$ est nulle (section 3.1). Les travaux récents de M. McQuillan déterminent la nature de ces feuilletages ; ils sont quotient d'un feuilletage induit par un champ de vecteurs, ce qui nous permet de conclure. Si la croissance de $\varphi$ est quadratique, il existe une fibration elliptique invariante qui s'avère être le seul feuilletage $\varphi$-invariant (section 3.2). Enfin, si tous les éléments de $\text{Aut}(\mathcal{F})$ ont une croissance des degrés bornée, nous montrons que $\mathcal{F}$ est invariant par un champ de vecteurs holomorphe.

Pour démontrer le théorème 1.2, nous nous appuierons sur les propriétés du feuilletage $\mathcal{F}$ sans utiliser directement les résultats de M. McQuillan. Suivant [12], il existe deux classes d'applications birationnelles non conjuguées à des automorphismes : ou bien la croissance des degrés est linéaire et la transformation préserve une fibration rationnelle, ou bien la croissance des degrés est exponentielle (voir la section 2). Nous montrons tout d'abord que, si $\varphi$ préserve une fibration rationnelle, $\varphi$ ne peut préserver d'autre feuilletage (section 4). Lorsque la croissance des degrés est exponentielle, nous proposons une étude basée sur les arguments de [6]. Le feuilletage $\mathcal{F}$ doit en effet être de type Riccati pour deux fibrations rationnelles transverses et nous en déduisons que la situation est conjuguée à celle de l'exemple 1.3 (section 5.5).

1.5. **Plan de l'article.** La section 2 contient les préliminaires nécessaires dans la suite ; nous conseillons au lecteur de lire le passage concernant la classification des applications birationnelles. Les trois sections suivantes démontrent les résultats principaux de ce texte. À la section 6, nous décrivons succinctement deux autres approches permettant de démontrer les points-clef du théorème 1.2. La section 7 regroupe les corollaires qui nous paraissent les plus frappants.

Nous avons délibérément choisi de ne pas détailler la théorie de la classification des feuilletages holomorphes singuliers. Nous renvoyons à l'excellent monographe de M. Brunella [7]. Cet article n'aurait certainement pas vu le jour sans l'existence de ce très bel ouvrage.

1.6. **Remerciements.** Cet article doit beaucoup aux discussions fructueuses que nous avons eues avec M. Brunella, D. Cerveau, S. Lamy, L. G. Mendes, J. V. Pereira et F. Touzet. Nous les remercions pour leurs conseils et leur enthousiasme.

## 2. Généralités sur les applications biméromorphes et les feuilletages

2.1. **Applications biméromorphes et stabilité algébrique.** Une application *biméromorphe* $\varphi$ entre deux espaces analytiques complexes $X$ et $Y$ est définie par son graphe $\Gamma \subset X \times Y$, une sous-variété analytique de $X \times Y$ telle que les deux projections $\pi_1 : \Gamma \to X$ et $\pi_2 : \Gamma \to Y$ soient des modifications propres (voir [13]). Lorsque $X$ et $Y$ sont des variétés projectives, $\varphi$ est dite *birationnelle*. Il sera souvent plus agréable de considérer une désingularisation $p : \widehat{\Gamma} \to \Gamma$ de $\Gamma$ ; nous noterons alors $p_1 = \pi_1 \circ p$, et $p_2 = \pi_2 \circ p$.

En dimension 2, toute modification propre est la composition d'un nombre fini d'éclatements de points (voir par exemple [24]). Toute application biméromorphe $\varphi$ se factorise donc sous la forme d'une composition finie d'éclatements suivie d'une suite finie de contractions de courbes exceptionnelles. En particulier, l'ensemble critique de $\varphi$ est une union de courbes rationnelles, éventuellement singulières. Par définition, l'ensemble des points de $X$ qui sont éclatés est l'ensemble des points d'indétermination : ce sont les points qui sont envoyés par $\varphi$ sur des courbes.

Décrivons rapidement les principaux résultats de [10] et [12] concernant la dynamique des applications biméromorphes qui seront utilisés dans la suite. Fixons $\varphi : X \circlearrowleft$, une application biméromorphe d'une surface complexe compacte *kählérienne*. L'application $\varphi$ induit une action naturelle $\varphi^* : H^{1,1}(X) \circlearrowleft$ définie par composition $\varphi^* := p_{1*} \circ p_2^*$. D'après les travaux de M. Gromov (voir [21]), le rayon spectral de cette transformation linéaire majore l'entropie topologique de $\varphi$. Il est donc naturel d'étudier le comportement asymptotique de la suite $\{\|(\varphi^n)^*\|\}_{n \geq 0}$ (pour une norme $\|.\|$ quelconque) ; celui-ci ne change pas si l'on conjugue $\varphi$ à l'aide d'un changement de variables birationnel (voir [12]) et détermine donc un invariant birationnel important de $\varphi$.

La notion centrale permettant d'étudier cette suite a été introduite dans [14].

**Définition 2.1.** Une application (bi)-méromorphe $\varphi : X \circlearrowleft$ est dite algébriquement stable (AS en abrégé) si l'une des conditions équivalentes suivantes est satisfaite :

a. Pour toute courbe $V$, il ne peut exister d'entier $n \geq 0$ tel que $\varphi^n(V)$ soit inclus dans l'ensemble d'indétermination de $\varphi$.
b. Pour tout entier $n \geq 0$, les deux transformations linéaires $(\varphi^n)^*$ et $(\varphi^*)^n$ de $H^{1,1}(X)$ coïncident.

Lorsque $\varphi$ est algébriquement stable, l'étude de la suite $\{\|(\varphi^n)^*\|\}_{n \geq 0}$ se ramène donc à l'étude du spectre de l'application linéaire $\varphi^*$ (point b. de la définition). Rendre une application algébriquement stable c'est, en quelque sorte, réduire les singularités qui apparaissent quand $\varphi$ est itérée. En ce sens, la stabilité algébrique est analogue à la notion de singularité réduite pour les feuilletages. Nous pouvons donc considérer le théorème suivant comme l'analogue du théorème de Seidenberg que nous présenterons au paragraphe 2.3.

**Théorème 2.1** ([12]). *Si $\varphi : X \circlearrowleft$ est une transformation biméromorphe d'une surface complexe compacte, il existe un morphisme biméromorphe $\pi : Y \to X$ tel que le relevé $\pi^{-1} \circ \varphi \circ \pi$ soit une transformation algébriquement stable de $Y$.*

**2.2. Classification suivant la croissance des degrés.** Il est possible de classer les transformations biméromorphes suivant la croissance des degrés, i.e.

suivant le comportement asymptotique de $\|(\varphi^*)^n\|$. Pour présenter cette classification, nous supposerons dans toute la suite de la discussion que $\varphi$ est algébriquement stable.

Supposons tout d'abord que $\varphi$ est un automorphisme ; $\varphi^*$ préserve alors la forme d'intersection sur $H^{1,1}(X)$. Puisque $X$ est kählérienne, celle-ci est non-dégénérée et possède une et une seule valeur propre positive (c'est le théorème de l'indice de Hodge). Si le spectre de $\varphi^*$ n'est pas inclus dans le disque unité de $\mathbb{C}$, il existe alors une unique valeur propre $\lambda$ de module strictement plus grand que 1 : c'est une valeur propre réelle simple.

Dans le cas général, une application biméromorphe ne préserve pas la forme d'intersection mais la dilate (voir [12]). Ceci permet d'obtenir des résultats analogues ; en particulier, la suite $\|(\varphi^*)^n\|$ ne peut se comporter que de quatre façons distinctes qui, suivant les cas, conduisent à la liste exhaustive suivante (pour les automorphismes voir [20] et [10]):

1.– Lorsque $\|(\varphi^n)^*\| = O(1)$, il existe un entier $n$ strictement positif tel que $\varphi^n$ est le flot au temps 1 d'un champs de vecteur (éventuellement trivial). Tout automorphisme d'une surface rationnelle minimale, ou d'une surface de type général est de ce type.

2.– Si $\|(\varphi^n)^*\| \sim C^{ste} \cdot n$, l'application $\varphi$ est birationnellement conjuguée à une application préservant une fibration rationnelle. En outre, $\varphi$ n'est pas birationnellement conjuguée à un automorphisme. L'application $(z, w) \to (z, zw)$ dans $\mathbb{P}^1 \times \mathbb{P}^1$ est un exemple de ce type.

3.– Lorsque $\|(\varphi^n)^*\| \sim C^{ste} \cdot n^2$, l'application $\varphi$ est birationnellement conjuguée à un automorphisme préservant une unique fibration elliptique. L'exemple le plus simple est donné par $(z, w) \to (z, z + w)$ dans le produit de deux courbes elliptiques identiques $X = \mathbb{C}/\Lambda \times \mathbb{C}/\Lambda$.

4.– Si $\|(\varphi^n)^*\| \sim \lambda^n$ avec $\lambda > 1$, le spectre de $\varphi^*$ contient une unique valeur propre hors du disque unité, à savoir $\lambda$, et celle-ci est simple. Il existe en outre un vecteur propre $\theta_+ \in \mathrm{H}^{1,1}(X, \mathbb{R})$ qui est numériquement effectif, ce qui signifie qu'il appartient à l'adhérence des classes engendrées par les formes de Kähler.

4.a.– Lorque $\theta_+ \cdot \theta_+ = 0$, l'application $\varphi$ est birationnellement conjuguée à un automorphisme. Dans ce cas, $X$ est une surface rationnelle, un tore, une surface K3 ou une surface de Enriques. L'exemple le plus simple est fourni par les automorphismes linéaires Anosov sur les tores.

4.b.– Lorsque $\theta_+ \cdot \theta_+ > 0$, l'application $\varphi$ n'est pas birationnellement conjuguée à un automorphisme. Dans ce cas, la surface $X$ est rationnelle. Les transformations de Hénon appartiennent à cette dernière classe.

Concluons cette section par le lemme suivant qui sera utilisé à plusieurs reprises.

**Lemme 2.2.** *Si $\varphi$ préserve deux fibrations distinctes, l'action en cohomologie vérifie $(\varphi^n)^* = Id$ pour un entier $n \geq 1$.*

2.3. **Quelques rappels sur les singularités des feuilletages.** Nous rappelons dans ce paragraphe les principaux résultats sur les singularités des feuilletages qui nous seront utiles pour la suite. L'exposition est, de ce fait, taillée sur mesure pour nos besoins. Nous renvoyons à [7] pour un texte précis et détaillé sur ce sujet.

Soit $\mathcal{F}$ un feuilletage holomorphe donné au voisinage d'un point $p$ par une 1-forme holomorphe $\omega$. En coordonnées locales, nous noterons
$$\omega = a(x,y)dx + b(x,y)dy$$
et $p$ sera identifié à l'origine de $\mathbb{C}^2$ et nous pouvons toujours supposer que les courbes $a^{-1}(0)$ et $b^{-1}(0)$ ne possèdent aucune composante commune. Le point $p$ est une singularité du feuilletage lorsque $a(0) = b(0) = 0$.

Si l'origine est éclatée, le feuilletage se relève en un nouveau feuilletage. Ce dernier peut alors présenter de nombreuses singularités le long du diviseur exceptionnel, mais celles-ci sont "plus simples" que la singularité initiale. Pour préciser cela, il convient d'introduire quelques définitions.

Une *séparatrice* de $\mathcal{F}$ en un point $p$ est un germe de courbe analytique passant par $p$ et tangent à $\mathcal{F}$, ce qui signifie que la forme $\omega$ définissant $\mathcal{F}$ est identiquement nulle le long du germe. Par définition, un feuilletage présente une *singularité dicritique* en un point $p$ s'il existe une infinité de séparatices passant par $p$. Lorsqu'une feuille de $\mathcal{F}$ contient une séparatrice, on dit également que cette feuille est une séparatrice.

Le point $p$ est une singularité *non-dégénérée* de $\mathcal{F}$ si la partie linéaire de $\omega$ possède deux valeurs propres non nulles dont le rapport n'est pas un nombre rationnel positif ; autrement dit, après changement de variables linéaire,
$$\omega = xdy + \lambda y dx + O(2), \tag{2.1}$$
avec $\lambda$ dans $\mathbb{C} \setminus \mathbb{Q}_+$. Le feuilletage possède exactement deux séparatrices passant en $p$, chacune d'entre elles étant tangente à l'un des axes de coordonnées.

Lorsque la partie linéaire de $\omega$ possède une et une seule valeur propre non nulle, l'origine est une singularité de type *selle-noeud*. Dans ce cas, il existe un changement de coordonnées holomorphes tel que
$$\omega = [x(1 + \nu y^k) + yF(x,y)]dy - y^{k+1}dx \tag{2.2}$$
où $\nu$ est un nombre complexe et $F$ est une fonction holomorphe qui s'annule à l'ordre $k$ en $p$ (forme normale de Dulac). Le feuilletage $\mathcal{F}$ possède toujours une séparatrice tangente à l'axe $\{y = 0\}$, appelée *séparatrice forte*. Il existe une séparatrice formelle tangente à $\{x = 0\}$, dite *séparatrice faible*, mais en général celle-ci n'est pas convergente.

Par définition, une singularité est *réduite* si elle est non-dégénérée ou de type selle-noeud. Une singularité dicritique n'est donc pas réduite. Un feuilletage est réduit si toutes ses singularités le sont.

**Théorème 2.3** (Seidenberg). *Soit $\mathcal{F}$ un feuilletage holomorphe de dimension 1 sur une surface complexe compacte $X$. Il existe alors un morphisme birationnel $\pi : Y \to X$ tel que $\pi^*\mathcal{F}$ soit un feuilletage réduit de $Y$.*

Lorsqu'on éclate une singularité réduite, il y a deux singularités le long du diviseur exceptionnel situées au niveau des transformées strictes des séparatrices. Un calcul direct montre que ces deux singularités sont réduites (voir aussi le paragraphe 4.1). De même, si l'on éclate un point lisse de $\mathcal{F}$, il apparaît une unique singularité le long du diviseur exceptionnel et celle-ci est réduite. Lorsqu'on effectue des éclatements, un feuilletage réduit reste donc réduit. Le théorème de Seidenberg peut donc être relu de la manière suivante : si l'on éclate toutes les singularités d'un feuilletage et que l'on itère ce procédé, on obtient un feuilletage réduit au bout d'un nombre fini d'étapes.

Nous aurons aussi besoin du théorème de Camacho-Sad. Nous ne l'énonçons que dans le cadre particulier où $\mathcal{F}$ ne possède que des singularités réduites. Supposons donc que $\mathcal{F}$ possède une singularité réduite en $p$, et considérons $C$ une séparatrice (nécessairement lisse) de $\mathcal{F}$ en $p$. L'indice de Camacho-Sad, noté $\mathbf{CS}(C, \mathcal{F}, p)$, peut être défini comme suit.

– Lorsque $p$ est une singularité non-dégénérée, donnée par la forme $\omega$ de la formule 2.1, $\mathbf{CS}(C, \mathcal{F}, 0) = -1/\lambda$ si $C$ est tangente à l'axe des $y$ et $\mathbf{CS}(C, \mathcal{F}, 0) = -\lambda$ si $C$ est tangente à l'axe des $x$.

– Lorsque $p$ est de type selle-noeud, donné par la forme $\omega$ de l'équation 2.2, $\mathbf{CS}(C, \mathcal{F}, 0) = 0$ si $C$ est la séparatrice forte, et $\mathbf{CS}(C, \mathcal{F}, 0) = \nu$ sinon.

Notons que, lorsqu'un éclatement de centre $p$ est réalisé, l'indice de Camacho-Sad des transformées strictes des séparatrices de $p$ chute d'une unité au point d'intersection avec le diviseur exceptionnel.

**Théorème 2.4** (Camacho-Sad). *Soit $\mathcal{F}$ un feuilletage holomorphe et $C$ une courbe compacte lisse invariante par $\mathcal{F}$. Si les singularités de $\mathcal{F}$ le long de $C$ sont réduites, alors*
$$[C]^2 = \sum_{Sing(\mathcal{F})} \mathbf{CS}(C, \mathcal{F}, p) .$$

Nous avons vu au paragraphe 1.1 comment associer à $\mathcal{F}$ son fibré cotangent $T_\mathcal{F}^*$. Le dual $T_\mathcal{F}$ est le fibré tangent du feuilletage. Lorsque $V$ est une courbe compacte, il est possible de calculer l'intersection $[T_\mathcal{F}] \cdot [V]$ à l'aide de formules

d'indices :

$$[T_\mathcal{F}] \cdot [V] = \chi(V) - Z(\mathcal{F}, V) \text{ si } V \text{ est } \mathcal{F}\text{-invariante}, \quad (2.3)$$
$$[T_\mathcal{F}] \cdot [V] = V^2 - \text{Tang}(\mathcal{F}, V) \text{ sinon}. \quad (2.4)$$

Le nombre $\chi(V)$ désigne la caractéristique d'Euler-Poincaré de $V$. Le nombre de points de tangence, compté avec multiplicité, est noté $\text{Tang}(\mathcal{F}, V)$. C'est un nombre entier positif.

De même, $Z(\mathcal{F}, V)$ est un entier qui est la somme des $Z(\mathcal{F}, V, p)$ où $p$ parcourt l'ensemble des singularités de $\mathcal{F}$ sur $V$ ; lorsque $p$ est non-dégénéré $Z(\mathcal{F}, V, p)$ vaut 1, lorsque $p$ est selle-noeud, $Z(\mathcal{F}, V, p)$ vaut 1 si $V$ est la séparatrice forte et $k$ sinon (voir l'équation (2.2)).

## 3. Le cas des automorphismes

Dans cette partie, nous montrons le théorème 1.1. Dans les deux premiers paragraphes, nous fixons un automorphisme d'ordre infini $\varphi$ qui préserve un feuilletage $\mathcal{F}$. Nous supposons d'abord que la croissance des degrés de $\varphi$ est exponentielle, puis qu'elle est quadratique. Le dernier paragraphe traite le cas des automorphismes isotopes à l'identité et termine la preuve du théorème 1.1.

### 3.1. Automorphismes dont l'entropie est positive.
Le but de ce paragraphe est de montrer le théorème suivant.

**Théorème 3.1.** *Soit $\varphi : X \to X$ un automorphisme d'une surface complexe projective dont la croissance des degrés est exponentielle. Si $\varphi$ préserve un feuilletage holomorphe $\mathcal{F}$, il existe alors un tore $T$ et un groupe fini $G$ d'automorphismes de $T$ tels que*

(i) *$X$ est isomorphe à la désingularisée de $T/G$, éventuellement éclatée ;*
(ii) *$\varphi$ se relève à $T$ en un automorphisme affine de type Anosov ;*
(iii) *le feuilletage $\mathcal{F}$ se relève en le feuilletage linéaire stable (ou instable) de ce difféomorphisme Anosov.*

**Remarque 3.1.** Dire que la croissance des degrés est exponentielle signifie que le rayon spectral de la transformation linéaire $\varphi^* : H^{1,1}(X) \circlearrowleft$ est strictement plus grand que 1. D'après les travaux de M. Gromov et Y. Yomdin ([21], [31]), cette propriété est équivalente à la stricte positivité de l'entropie topologique de $\varphi$. L'hypothèse de projectivité n'est pas nécessaire. Elle sera seulement utilisée pour appliquer le théorème de Miyaoka sur la positivité du fibré cotangent de $\mathcal{F}$.

**Corollaire 3.2.** *Soit $\varphi : X \to X$ un automorphisme d'une surface projective dont l'entropie topologique est strictement positive. Si $\varphi$ préserve un feuilletage holomorphe, $\varphi$ en préserve deux.*

**Exemple 3.1.** Soit $T$ un tore complexe de dimension 2 admettant des automorphismes linéaires de type Anosov (nous renvoyons à [19] pour la classification de ces tores). Chacun de ces automorphismes commute avec l'involution $\sigma(x,y) = (-x,-y)$ et passe donc au quotient sur $T/\langle\sigma\rangle$. Après résolution des singularités, nous obtenons ainsi une surface K3 munie d'automorphismes d'entropie positive préservant chacun deux feuilletages holomorphes.

**Exemple 3.2.** Soit $\xi$ une racine de l'unité d'ordre 3, 4 ou 6 et $E$ la courbe elliptique $\mathbb{C}/\mathbb{Z}[\xi]$. La représentation linéaire standard de $\mathrm{SL}(2,\mathbb{Z}[\xi])$ sur $\mathbb{C}^2$ induit un plongement naturel de $\mathrm{SL}(2,\mathbb{Z}[\xi])$ dans les automorphismes de la variété abélienne $T = E \times E$. Tous les automorphismes ainsi obtenus commutent avec l'homothétie $h : T \to T$ donnée par $h(x,y) = (\xi x, \xi y)$. La surface $S$ obtenue en désingularisant le quotient $T/\langle h \rangle$ est une surface rationnelle. Tous les automorphismes issus de $\mathrm{SL}(2,\mathbb{Z}[\xi])$ passent au quotient et, en dehors de ceux provenant de matrices unipotentes, déterminent des automorphismes de $S$ préservant deux feuilletages.

Le Théorème 3.1 implique le corollaire suivant, dont nous donnons une preuve détaillée dans l'appendice A.

**Corollaire 3.3.** *Soit $\varphi : X \to X$ un automorphisme d'entropie positive sur une surface rationnelle. Si $\varphi$ préserve un feuilletage holomorphe, la surface est biholomorphe, après un nombre fini de contractions $\varphi$-équivariantes, à l'une des surfaces de l'exemple 3.2.*

Le reste de cette section est occupé par la preuve du Théorème 3.1. Avant de considérer le cas général, traitons tout d'abord le cas des feuilletages engendrés par un champ de vecteurs.

**Proposition 3.4.** *Soit $v$ un champ de vecteurs holomorphe sur une surface complexe compacte, et $\mathcal{F}$ le feuilletage associé. Si $\mathcal{F}$ est invariant par un automorphisme d'entropie positive, alors $X$ est un tore et $v$ est donc linéaire.*

*Démonstration.* Toutes les courbes exceptionnelles de première espèce tracées sur $X$ ont une auto-intersection strictement négative et sont donc invariantes par $v$. S'il en existe une infinité, un théorème de Jouanolou montre que le feuilletage induit par $v$ est une fibration rationnelle préservée par $\varphi$ (voir par exemple [18]). Ceci contredit l'hypothèse d'entropie positive. Il n'y a donc qu'un nombre fini de courbes exceptionnelles ; elles sont permutées par $\varphi$. Il s'ensuit que $v$, $\mathcal{F}$ et $\varphi$ descendent sur un modèle minimal $X'$ de $X$ ; nous les noterons encore $\varphi$, $v$ et $\mathcal{F}$. Puisque la croissance des degrés est invariante par changement de variable birationnel (section 2), l'automorphisme induit par $\varphi$ sur $X'$ reste d'entropie positive.

Parmi les surfaces minimales, seuls les tores, les surfaces K3 et les surfaces d'Enriques possèdent des automorphimes d'entropie positive (voir [9]). Les surfaces K3 et les surfaces d'Enriques ne possèdent aucun champ de vecteurs. Nous avons donc montré que $X'$ est un tore, que $\varphi$ est un difféomorphisme affine de type Anosov et que $\mathcal{F}$ est linéaire. Puisque $\mathcal{F}$ est linéaire, $v$ n'a pas de zéros et ne se relève à aucun éclatement de $X'$. Ceci montre que $X$ et $X'$ coïncident. La proposition est démontrée. □

*Preuve du théorème 3.1.* Soit $\varphi : X \to X$ un automorphisme d'entropie positive préservant un feuilletage holomorphe $\mathcal{F}$. Nous allons montrer que, quitte à relever $\varphi$ et $\mathcal{F}$ sur un revêtement d'ordre fini de $X$, le feuilletage $\mathcal{F}$ est engendré par un champ de vecteurs. La conclusion découlera alors de la proposition précédente. La preuve est une application directe de la classification des feuilletages de dimension de Kodaira nulle. Nous en rappelons les ingrédients essentiels pour la commodité du lecteur.

Les singularités de $\mathcal{F}$ sont en nombre fini et sont permutées par $\varphi$. Le théorème de Seidenberg permet donc de supposer que $\mathcal{F}$ est réduit (section 2). Les courbes exceptionnelles de première espèce qui sont invariantes par $\mathcal{F}$ sont permutées par $\varphi$. Nous pouvons donc contracter celles qui ne créent pas de singularités non-réduites par contraction : le feuilletage $\mathcal{F}$ devient *relativement minimal* dans la terminologie de Brunella-McQuillan.

Comme $\mathcal{F}$ est $\varphi$-invariant et $\varphi$ est un automorphisme,

$$\varphi^* T_\mathcal{F}^* = T_\mathcal{F}^*, \tag{3.1}$$

ce qui va nous permettre de décrire les propriétés numériques de $T_\mathcal{F}^*$. Puisque $X$ est projective, le théorème de Miyaoka s'applique : ou bien $\mathcal{F}$ est une fibration rationnelle, ou bien $T_\mathcal{F}^*$ est pseudo-effectif, c'est-à-dire que $T_\mathcal{F}^* \cdot C \geq 0$ pour tout diviseur ample $C$ (voir [7]). Puisque l'entropie topologique de $\varphi$ est strictement positive, nous sommes dans ce dernier cas. D'après [15], il existe une décomposition de Zariski

$$c_1(T_\mathcal{F}^*) = [P] + [N], \tag{3.2}$$

où $P$ est un $\mathbb{Q}$-diviseur numériquement effectif, $N = \sum_i a_i D_i$ est un $\mathbb{Q}^+$-diviseur contractible, c'est-à-dire que la matrice d'intersection $(D_i \cdot D_j)$ est définie négative, et $P \cdot D_i = 0$ pour tout $i$. Suivant [7], nous savons que

- chaque composante connexe du support de $N$ est une chaîne de courbes rationnelles s'intersectant transversalement,
- les coefficients $a_i$ sont strictement plus petits que 1.

Le diviseur $N$ est uniquement déterminé par $\mathcal{F}$, donc $\varphi^* N = N$, et par suite $\varphi^*[P] = [P]$ ; par conséquent, $[P]$ est orthogonal à $\mathbb{R}\theta_+$, l'unique direction dilatée par $\varphi^*$. Comme $[P]^2 \geq 0$ et $\theta_+^2 \geq 0$, le théorème de l'indice de Hodge montre que $[P]$ est proportionnel à $\theta_+$ et donc que $[P] = 0$. Dans la terminologie de Brunella-McQuillan, nous avons montré que la dimension de Kodaira numérique de $\mathcal{F}$ est nulle.

Comme $X$ est une surface rationnelle, une surface de Enriques ou une surface K3 (si $X$ est un tore, il n'y a rien à démontrer), nous savons que le premier nombre de Betti de $X$ est nul. Le fibré $T_{\mathcal{F}}^{*\otimes n}$ admet donc une section non nulle si $n$ désigne le plus petit commun multiple des dénominateurs des $a_i$. Notons $E_1$ l'espace total du fibré $T_{\mathcal{F}}^*$ et $E_n$ celui de $T_{\mathcal{F}}^{*\otimes n}$. Soient $s : X \to E_n$ la section de $T_{\mathcal{F}}^{*\otimes n}$, $p_n : E_1 \to E_n$ le morphisme d'élévation à la puissance $n$ et $\Phi^{\otimes n} : E_n \to E_n$ l'automorphisme induit par le relevé $\Phi$ de $\varphi$ à $E_1$ ; en choisissant $\Phi$ convenablement nous pouvons supposer que $\Phi^{\otimes n} \circ s = s \circ \varphi$.

Notons $Y_0$ l'image réciproque de $s(X)$ par $p_n$ et $Y$ la désingularisation minimale de $Y_0$. Puisque $\Phi^{\otimes n} \circ s = s \circ \varphi$ et $p_n \circ \Phi = \Phi^{\otimes n} \circ p_n$, l'automorphisme $\Phi$ induit par restriction un automorphisme de $Y_0$. Celui-ci se relève canoniquement en un automorphisme $\widehat{\Phi}$ de $Y$ [24].

Soit $\pi : Y \to X$ la composée du morphisme de désingularisation $Y \to Y_0$ et de la projection $Y_0 \to X$. L'image réciproque de $\mathcal{F}$ par $\pi$ est un feuilletage holomorphe $\mathcal{G}$ dont le fibré canonique coïncide avec $\pi^*(T_{\mathcal{F}}^*)$ (voir [7]). Par construction, ce fibré admet une section holomorphe non nulle $\widehat{s}$.

A priori, le feuilletage induit sur $Y$ n'est pas relativement minimal et l'on peut contracter des courbes exceptionnelles de première espèce tangentes à $\mathcal{G}$ (donc permutées par $\widehat{\Phi}$) apparues lors de la désingularisation. Notons $\epsilon : Y \to A$ le morphisme birationnel correspondant à cette série de contractions. L'automorphisme $\widehat{\Phi}$ passe au quotient en un nouvel automorphisme préservant le feuilletage $\epsilon_*\mathcal{G}$ ; puisque le fibré canonique de ce feuilletage admet une section holomorphe, à savoir $\epsilon_*(\widehat{s})$, ce fibré est trivial car son diviseur des zéros est du type précédent $\sum a_i D_i$ avec $a_i < 1$ ; le morphisme $\epsilon$ contracte donc le lieu des zéros de $\hat{s}$. Le feuilletage $\epsilon_*\mathcal{G}$ est donc induit par un champ de vecteurs global. La proposition précédente montre que $A$ est un tore. $\square$

3.2. **Automorphismes préservant une fibration elliptique.** Dans ce paragraphe, nous nous intéressons au cas des automorphismes qui ont une entropie nulle et dont aucun itéré n'est le flot d'un champ de vecteurs. Nous avons vu à la partie 2 que ceci est équivalent à une croissance quadratique des degrés : $\|\varphi^{n*}\| \sim C^{ste} \cdot n^2$.

**Théorème 3.5.** *Soit $\varphi : X \to X$ un automorphisme d'une surface complexe compacte kählérienne. Si $\|\varphi^{n*}\| \sim C \cdot n^2$, alors $\varphi$ préserve une fibration elliptique et cette fibration est le seul feuilletage holomorphe de $X$ invariant par $\varphi$.*

D'après le paragraphe 2.2, nous savons que $\varphi$ préserve une unique fibration elliptique $\pi : X \to B$. Notons $\overline{\varphi} : B \to B$ l'action induite par $\varphi$ dans la base.

Lorsque $\overline{\varphi}$ est d'ordre fini, $\varphi$ agit par translation dans chaque fibre et la condition de croissance $\|\varphi^{n*}\| \sim C^{ste} \cdot n^2$ force le paramètre de translation à dépendre non-trivialement de la base. L'automorphisme $\varphi$ crée donc un "twist", ce qui force tout feuilletage $\varphi$-invariant à être parallèle aux fibres. Pour mettre en œuvre ces idées heuristiques, nous établirons la

**Proposition 3.6.** *Soit $\varphi$ un automorphisme d'une surface complexe compacte kählérienne $X$ vérifiant $\|\varphi^{n*}\| \sim C^{ste} \cdot n^2$. Alors $\varphi$ préserve une unique fibration elliptique $\pi : X \to B$ et, si le modèle minimal de $X$ n'est pas un tore, l'action de $\varphi$ dans la base $B$ est d'ordre fini.*

*Démonstration.* Il s'agit de montrer la seconde assertion, celle relative à l'action de $\overline{\varphi}$. Supposons donc que $X$ n'est pas un tore et que $\overline{\varphi}$ est d'ordre infini. Cette dernière hypothèse impose à $B$ d'être elliptique ou rationnelle ; les valeurs critiques de $\pi$ étant permutées par $\overline{\varphi}$, il n'y en a pas si $B$ est elliptique et il y en a au plus 2 si $B$ est rationnelle. Quitte à remplacer $\varphi$ par l'un de ses itérés, nous pouvons supposer que toutes les fibres multiples ou singulières sont fixes par $\varphi$.

Commençons par le cas où $B$ est elliptique. Puisqu'aucune fibre n'est singulière, l'invariant modulaire des fibres de $\pi$ détermine une fonction holomorphe sur $B$ ; il s'agit donc d'une constante et $X$ est une fibration localement triviale. Autrement dit, $X$ est une suspension d'une courbe elliptique $E$ au-dessus de $B$. Puisque $\mathrm{Aut}(E)$ est virtuellement abélien, $X$ est le quotient d'un tore par un groupe fini de transformations sans point fixe. Si ce quotient n'est pas un tore, $X$ est une surface bi-elliptique et possède une seconde fibration invariante par $\mathrm{Aut}(X)$ [9]. Ceci contredit la condition de croissance $\|\varphi^{n*}\| \sim C^{ste} \cdot n^2$.

Supposons maintenant que $B$ est isomorphe à la droite projective $\mathbb{P}^1$. L'action de $\varphi$ dans la base est une translation $\overline{\varphi}(z) = z+1$ ou une homothétie $\overline{\varphi}(z) = \alpha z$ d'ordre infini. Puisque l'invariant modulaire des fibres est une fonction méromorphe invariante sous l'action de $\overline{\varphi}$, il doit être constant ; autrement dit, la fibration elliptique est isotriviale. En effectuant un changement de base $z \mapsto z^n$ pour un entier $n \geq 2$ convenable, les fibres singulières sont alors remplacées par des fibres non-multiples de type $I_b$, $b \geq 0$ (terminologie de Kodaira [2]). Lorsque $b$ est strictement positif, l'invariant modulaire possède un pôle d'ordre $b$ et la fibration ne peut être isotriviale. Le produit fibré

$$Y = \{(x,z) \in X \times \mathbb{P}^1 \mid \pi(x) = z^n\}$$

obtenu par le changement de base est donc une fibration elliptique localement triviale sur $\mathbb{P}^1$, donc triviale. Autrement dit, $Y$ est isomorphe à $\mathbb{P}^1 \times E$ pour une courbe elliptique $E$.

Notons $p : Y \to X$ la projection induite par le changement de base $z \mapsto z^n$. L'image par $p$ d'une courbe rationnelle $\mathbb{P}^1 \times \{m\}$, avec $m$ générique, est une courbe rationnelle lisse d'auto-intersection positive ou nulle. Donc, soit $X$ est rationnelle, soit $X$ est réglée à base elliptique. Le second cas est exclu par le Lemme 2.2, car $\varphi$ préserverait une fibration rationnelle et une fibration elliptique.

La surface $X$ est donc rationnelle. Le Lemme 3.7 ci-dessous montre que la fibration $\pi$ possède au moins deux fibres singulières. Celles-ci sont préservées par $\varphi$, donc $\overline{\varphi}(z) = \alpha z$ pour $\alpha \in \mathbb{C}^*$. Nous pouvons alors relever $\overline{\varphi}$ par le changement de base $z \mapsto z^n$ (prendre $\overline{\psi} = \beta z$ avec $\beta^n = \alpha$) : $\varphi$ se relève également en un automorphisme $\psi$ du produit $Y \cong \mathbb{P}^1 \times E$. Un tel automorphisme préserve automatiquement la fibration rationnelle. La suite $\|\psi^{n*}\|$ est donc bornée (lemme 2.2), ce qui est incompatible avec la croissance quadratique de la suite $\|\varphi^{n*}\|$. La proposition est donc démontrée. $\square$

**Lemme 3.7.** *Soit $\pi : X \to \mathbb{P}^1$ une fibration elliptique sur une surface rationnelle. Alors $\pi$ possède au moins deux fibres singulières.*

*Démonstration.* La preuve s'appuie sur un calcul de caractéristique d'Euler. Nous pouvons supposer qu'aucune fibre ne possède de courbes rationnelles de première espèce. Nous noterons $[F]$ la classe d'une fibre régulière. Si $X_z$ est une fibre singulière, $r_z$ désignera le nombre de ses composantes irréductibles et $\chi(X_z)$ désignera sa caractéristique d'Euler.

À partir de la table des fibres singulières de Kodaira (voir [2] p.150), il est facile de vérifier l'inégalité $\chi(X_z) \leq r_z + 1$. Comme les fibres lisses de la fibration ont une caractéristique d'Euler nulle, nous obtenons

$$\chi(X) = \sum_{X_z \text{ singulière}} \chi(X_z) \leq \sum_{X_z \text{ singulière}} (r_z + 1). \tag{3.3}$$

Par ailleurs, chaque composante irréductible $D_i$ d'une fibre singulière engendre une classe dans l'espace orthogonal $[F]^\perp$ de $[F]$ (pour la forme d'intersection). Si $X_z = \sum D_i$, la seule relation linéaire liant les $[D_i]$ est donnée par $\sum [D_i] = [F]$. Ceci montre l'inégalité

$$\sum (r_z - 1) \leq h^{1,1}(X) - 1. \tag{3.4}$$

Pour une surface rationnelle, nous savons en outre que

$$\chi(X) = 2 + h^{1,1}(X). \tag{3.5}$$

Si $N$ désigne le nombre de fibres singulières, les deux estimations précédentes entrainent

$$N + \sum r_z \geq 2 + h^{1,1}(X) \tag{3.6}$$

$$\geq 3 - N + \sum r_z, \tag{3.7}$$

ce qui montre que $N$ est supérieur ou égal à 2. □

Déduisons maintenant le Théorème 3.5 de la proposition précédente.

*Démonstration du théorème 3.5.* Soit $\varphi : X \to X$ un automorphisme vérifiant la condition de croissance $\|\varphi^{n*}\| \sim C^{ste} \cdot n^2$ et préservant un feuilletage $\mathcal{F}$. Lorsque $X$ est un tore, ou un tore éclaté, $\varphi$ est une application affine, $\mathcal{F}$ est linéaire et le théorème se démontre aisément. Lorsque $X$ n'est pas un tore, d'après la proposition 3.6, il existe une fibration elliptique $\pi : X \to B$ et un itéré de $\varphi$, encore noté $\varphi$ dans la suite, qui préserve chaque fibre de $\pi$. Nous noterons $[F]$ la classe de Chern de cette fibration.

Soit $\mathcal{F}$ un feuilletage $\varphi$-invariant ne coïncidant pas avec la fibration. Le lieu de tangence de $\mathcal{F}$ avec la fibration est un diviseur effectif $\varphi$-invariant (éventuellement nul) que nous noterons $T$.

Si $T$ possède une composante transverse aux fibres, chaque application $\varphi_z := \varphi|_{\pi^{-1}\{z\}}$ possède un point périodique. Il s'ensuit que $\varphi^n = \mathrm{Id}$ pour un entier $n \geq 1$, ce qui est absurde. Donc $\mathcal{F}$ est transverse aux fibres génériques de $\pi$. En particulier, la fibration est isotriviale. Quitte à faire un changement de base, nous pouvons donc supposer que $\pi$ est une fibration localement triviale sans fibre singulière.

Soit $\Delta \subset B$ un petit disque de la base. Il existe une fonction holomorphe $t : \Delta \to \mathbb{C}$, et un nombre $\tau \in \mathbb{C}\setminus\mathbb{R}$ tels que $\pi^{-1}(\Delta)$ soit isomorphe à $\Delta \times (\mathbb{C}/\mathbb{Z}+\tau\mathbb{Z})$ et tel que $\varphi(z, w) = (z, w + t(z))$. La différentielle de l'itéré $n$-ième de $\varphi$ est donnée par

$$d\varphi^n(z, w) = \begin{bmatrix} 1 & 0 \\ n\frac{dt}{dz} & 1 \end{bmatrix}. \tag{3.8}$$

En particulier, lorsque $dt/dz \neq 0$, les images d'un vecteur tangent quelconque par les itérés $d\varphi^n$ convergent vers un vecteur tangent à la fibration. Comme $\mathcal{F}$ est distinct de la fibration, $t(z) = c \in \mathbb{C}$ est constant.

L'automorphisme $\varphi$ peut alors être déformé continûment en une famille d'automorphismes $\varphi_\varepsilon$ fixant chaque fibre, telle que le paramètre de translation dans les fibres soit donné par $c + \varepsilon$ pour $\varepsilon \in \mathbb{C}$ arbitraire. Donc $\varphi^* = \varphi^*_{-c} = \mathrm{Id}$ ce qui contredit la condition de croissance, et termine la preuve. □

### 3.3. Automorphismes isotopes à l'identité.

Les paragraphes précédents permettent de décrire les feuilletages qui sont invariants par un automorphisme dont aucun itéré n'est isotope à l'identité. Nous allons maintenant traiter le cas restant, celui des feuilletages invariants par un automorphisme d'ordre infini qui coïncide avec le flot au temps 1 d'un champ de vecteurs holomorphes.

**Proposition 3.8.** *Soit $\mathcal{F}$ un feuilletage holomorphe singulier d'une surface projective $X$. Si $\mathcal{F}$ est préservé par le flot $\Phi : \mathbb{C} \times X \to X$ d'un champ de vecteurs holomorphes $v$, il s'agit d'une des cinq situations suivantes :*

(i) *le champ $v$ est tangent à un fibré elliptique et $\mathcal{F}$ est un feuilletage tourbillonné ou la fibration elle-même ;*
(ii) *la variété $X$ est un tore et $v$ et $\mathcal{F}$ sont linéaires ;*
(iii) *le champ $v$ est tangent à une fibration rationnelle et $\mathcal{F}$ est un feuilletage de Riccati à monodromie abélienne, ou la fibration elle-même ;*
(iv) *la surface $X$ est une suspension de $\mathbb{P}^1$ au-dessus d'une courbe elliptique, $v$ s'obtient également par suspension et $\mathcal{F}$ aussi, sauf s'il coïncide avec la fibration rationnelle ;*
(v) *quitte à faire un changement de variable birationnel $v$ est un champ de vecteurs linéaire de $\mathbb{P}^1 \times \mathbb{P}^1$ et $\mathcal{F}$ est donné par un champ linéaire qui commute avec $v$.*

*Preuve.* Conservons les notations de l'énoncé en notant $\Phi_t$ le flot de $v$ au temps $t$. La liste des champs de vecteurs holomorphes des surfaces complexes compactes est connue. Le livre [7] contient les détails de cette classification lorsque la surface est projective ; seuls cinq cas apparaissent :

**1.-** $X$ est un fibré elliptique, éventuellement avec des fibres multiples mais sans fibres singulières, et $v$ est tangent aux fibres. Plus précisément, il existe un revêtement ramifié de $X$ qui est biholomorphe à la suspension d'une courbe elliptique $E$ au-dessus d'une surface de Riemann $B$ et $v$ provient d'un champ constant sur $E$. Dans ce cas, les orbites du flot $\Phi$ sont toutes infinies et le feuilletage $\mathcal{F}$ n'a donc aucune singularité. Le lieu de tangence entre $\mathcal{F}$ et la fibration elliptique étant un diviseur invariant par $v$, c'est une union de fibres et $\mathcal{F}$ est un feuilletage tourbillonné (voir [17] pour une définition précise) ou la fibration elle-même.

**2.-** $X$ est un tore, $v$ est linéaire et nous pouvons supposer que les orbites du flot de $v$ sont Zariski-denses, car sinon le cas (1.-) peut-être appliqué. Le feuilletage $\mathcal{F}$ ne possède alors ni singularité ni feuille compacte. Il est facile d'en déduire que $\mathcal{F}$ est donné par les orbites d'un autre champ de vecteurs linéaire (voir [17]).

**3.-** $X$ est une fibration rationnelle et $v$ est tangent aux fibres de cette fibration. Supposons que $\mathcal{F}$ ne coïncide pas avec la fibration et montrons alors que le diviseur des tangences entre $\mathcal{F}$ et la fibration est une union de fibre.

S'il existe une courbe $V$, transverse à la fibration rationnelle, le long de laquelle $\mathcal{F}$ est tangent à la fibration, cette courbe est contenue dans l'ensemble des zéros de $v$. Soit $p$ un point générique de $V$. La feuille locale $L_p$ de $\mathcal{F}$ passant par $p$ est lisse, elle est invariante par le flot $\Phi_t$ de $v$ (car $p$ est fixé par le flot) et elle n'est pas contenue dans une fibre. En particulier, elle intersecte chaque fibre voisine de $p$ en un nombre fini de points. Comme $\Phi_t(L_p) = L_p$, ces points d'intersection sont permutés par le flot du champ ; ils sont donc en fait fixés par le flot, ce qui montre que $L_p$ est contenue dans les zéros de $v$. Puisque ceci est valable pour un choix générique de $p$ le long de la courbe $V$, $v$ serait nul dans un voisinage de $V$, donc partout nul.

Nous avons donc montré que le lieu de tangence avec la fibration rationnelle est constitué de fibres. Le feuilletage $\mathcal{F}$ est donc un feuilletage de Riccati. La monodromie est abélienne car, dans chaque fibre, elle commute avec le flot.

**4.-** $X$ est une suspension de $\mathbb{P}^1$ au-dessus d'une courbe elliptique $E$ et $v$ aussi. Autrement dit, $X$ est le fibré en droites projectives associé à une représentation $\rho : \pi_1(E) \to \mathrm{PGL}(2,\mathbb{C})$ ; le champ de vecteurs $v$ est obtenu par suspension : il se projette sur la courbe elliptique en un champ constant et non nul $v_E$. Il s'ensuit que $\mathcal{F}$ n'a pas de singularités. Nous supposerons de plus que les orbites de $v$ ne sont pas des courbes elliptiques, car sinon le cas (1.-) peut être appliqué. Si une des feuilles de $\mathcal{F}$ est une fibre rationnelle, $\mathcal{F}$ coïncide avec la fibration. Hormis ce cas, deux possibilités peuvent apparaître suivant la structure du lieu de tangence entre $\mathcal{F}$ et la fibration rationnelle :

a. $\mathcal{F}$ est partout transverse à la fibration : dans ce cas il est donné par un autre champ de vecteurs, par exemple la projection de $v$ sur $\mathcal{F}$ parallèlement à la fibration rationnelle.
b. $\mathcal{F}$ est tangent à la fibration rationnelle le long de certaines courbes.

Montrons que ce dernier cas peut être exclu. Pour cela, notons $\Lambda \subset \mathbb{C}$ le réseau des paramètres $t$ pour lesquels $\Phi_t$ fixe chaque fibre de la fibration rationnelle. Soit $p$ un point générique du diviseur de tangence entre $\mathcal{F}$ et la fibration, et $L_p$ la feuille locale de $\mathcal{F}$ en $p$. Le point $p$, et donc la feuille $L_p$, sont fixes sous l'action de $\{\Phi_t\}_{t \in \Lambda}$. Soit $F$ une fibre rationnelle voisine de la fibre passant par $p$. Les transformations $\Phi_t$, $t \in \Lambda$, sont des homographies de $F$ qui fixent les points de tangence de $\mathcal{F}$ avec la fibration et les points d'intersection de $L_p$ avec $F$ : ces homographies sont donc toutes périodiques, ceci quelque soit $t$ dans $\Lambda$. Par conséquent, la monodromie de $v$ est finie et les orbites de $v$ sont des courbes elliptiques. Ceci contredit nos hypothèses.

**5.-** Quitte à faire un changement de variables birationnel, $X$ est isomorphe à $\mathbb{P}^1 \times \mathbb{P}^1$ et $v = v_1 + v_2$ est linéaire.

Un champ de vecteurs sur la sphère de Riemann préserve toujours une 1-forme différentielle méromorphe fermée : si le flot est du type $(t, z) \mapsto \lambda^t z$, il suffit de prendre la forme $dz/z$, s'il est du type $(t, z) \mapsto z/(1 + tz)$, la forme $dz/z^2$ convient. Nous pouvons donc fixer deux formes méromorphes fermées $\Omega_1$ et $\Omega_2$ qui sont invariantes par le flot $\Phi_t$ associé au champ $v_1 + v_2$ et dont les noyaux sont respectivement tangents aux fibrations verticales et horizontales de $X = \mathbb{P}^1 \times \mathbb{P}^1$.

Soit $\Omega$ une 1-forme méromorphe définissant $\mathcal{F}$. Il existe alors deux fonctions méromorphes $a$ et $b$ telles que $\Omega = a\Omega_1 + b\Omega_2$. L'invariance de $\mathcal{F}$ et des $\Omega_i$ ($i=1,2$) par l'automorphisme $\Phi_t$ entraîne :

$$\frac{a}{b}(z) = \frac{a}{b} \circ \Phi_t(z), \quad \forall z \in \mathbb{P}^1. \tag{3.9}$$

Si $v$ n'est pas tangent à un pinceau de courbes rationnelles, cette équation montre que $a/b$ est une constante, ce qui signifie que $\mathcal{F}$ est un feuilletage linéaire associé à un champ de vecteurs qui commute à $v$. Lorsque $v$ est tangent à un pinceau de courbes rationnelles, nous pouvons effectuer un nombre fini d'éclatements qui nous ramènent dans une situation analogue au cas 3. □

3.4. **Preuve du théorème 1.2.** Pour terminer l'étude des surfaces feuilletées $(X, \mathcal{F})$ qui possèdent une infinité de symétries biholomorphes, nous aurons besoin de la proposition suivante:

**Proposition 3.9.** *Soit $(X, \mathcal{F})$ une surface projective feuilletée dont le groupe de symétries biholomorphes $\mathrm{Aut}(\mathcal{F})$ est infini. Alors $\mathrm{Aut}(\mathcal{F})$ contient un élément d'ordre infini. Plus précisément, soit $\mathrm{Aut}(\mathcal{F})$ contient un sous-groupe à un paramètre, soit il contient un élément $\varphi$ dont l'action sur $H^{1,1}(X)$ vérifie $\|\varphi^{n*}\| \to \infty$.*

*Démonstration.* Considérons la représentation

$$\rho : \mathrm{Aut}(X) \longrightarrow \mathrm{GL}(H_2(X, \mathbb{Z}))$$

donnée par l'action des automorphismes de $X$ sur son deuxième groupe d'homologie. Si $\mathrm{Aut}(\mathcal{F})$ contient un élément $\varphi$ pour lequel $\rho(\varphi)$ est infini, la suite $\|\varphi^{n*}\|$ tend vers l'infini quand $n$ tend vers l'infini. Nous pouvons donc supposer que $\rho(\mathrm{Aut}(\mathcal{F}))$ est un sous-groupe de torsion de $\mathrm{GL}(H_2(X, \mathbb{Z}))$. Un tel groupe est nécessairement fini, ce qui montre que $\mathrm{Aut}(\mathcal{F})$ possède un sous-groupe d'indice fini dont les éléments agissent trivialement sur l'homologie de $X$. D'après [25], il existe donc un sous-groupe d'indice fini $\mathrm{Aut}_0(\mathcal{F})$ dans $\mathrm{Aut}(\mathcal{F})$ dont tous les éléments sont holomorphiquement isotopes à l'identité.

Notons $\mathrm{Aut}_0(X)$ la composante connexe de l'identité du groupe de Lie complexe $\mathrm{Aut}(X)$, de sorte que

$$\mathrm{Aut}_0(\mathcal{F}) = \mathrm{Aut}_0(X) \cap \mathrm{Aut}(\mathcal{F}). \tag{3.10}$$

Si $\Omega$ est une 1-forme différentielle rationnelle définissant le feuilletage $\mathcal{F}$, le groupe $\mathrm{Aut}_0(\mathcal{F})$ est le sous-groupe de Lie complexe fermé de $\mathrm{Aut}_0(X)$ déterminé par l'équation

$$\varphi^*\Omega \wedge \Omega = 0. \tag{3.11}$$

Rappelons que le groupe $\mathrm{Aut}_0(X)$ scinde en une suite exacte

$$1 \to L \to \mathrm{Aut}_0(X) \to A \to 0, \tag{3.12}$$

où $L$ est un groupe algébrique linéaire et $A$ est un sous-groupe compact du groupe des translations de la variété d'Albanese de $X$ (voir [25] ou [8]). Si $\mathrm{Aut}_0(\mathcal{F}) \cap L$ est infini, par exemple si $A$ est nul, l'équation (3.11) montre que $\mathrm{Aut}_0(\mathcal{F}) \cap L$ est un sous-groupe algébrique infini de $L$. Un tel groupe ne possède qu'un nombre fini de composantes connexes ; dans ce cas, $\mathrm{Aut}_0(\mathcal{F})$ contient donc un groupe à un paramètre. Un argument analogue s'applique si $L$ est trivial : $\mathrm{Aut}_0(\mathcal{F})$ est alors un sous-groupe fermé infini de $A$, donc il contient un sous-groupe à un paramètre.

Dans le cas restant, les dimensions de $A$ et de $L$ sont strictement positives ; la surface $X$ est alors une surface réglée au-dessus d'une courbe elliptique qui possède des champs de vecteurs du type 4.- (voir la démonstration de la proposition 3.8 ci-dessus). Puisque $\mathrm{Aut}_0(\mathcal{F})$ est infini et qu'il ne coupe $L$ qu'en un nombre fini de points, l'action induite sur la base comporte une infinité de translations. Le feuilletage $\mathcal{F}$ est donc régulier. En reprenant la preuve de la proposition 3.8, point 4.-, il est facile de voir que $\mathcal{F}$ est invariant par un groupe à un paramètre. □

*Preuve du théorème 1.2.* Soit $\mathcal{F}$ un feuilletage holomorphe singulier sur une surface projective dont le groupe de symétries $\mathrm{Aut}(\mathcal{F})$ est infini. Appliquons la proposition 3.9. Si $\mathrm{Aut}(\mathcal{F})$ possède un sous-groupe à un paramètre, $\mathcal{F}$ est invariant par un champ de vecteurs et nous pouvons appliquer la proposition 3.8. Sinon $\mathrm{Aut}(\mathcal{F})$ possède un élément d'ordre infini, dont aucun itéré n'est le flot d'un champ de vecteurs car $\|\varphi^{n*}\| \to \infty$. Si la croissance de $\|\varphi^{n*}\|$ est exponentielle, le théorème 3.1 peut être appliqué : la situation est celle de l'exemple 1.1. Sinon la croissance de $\|\varphi^{n*}\|$ est quadratique et le théorème 3.5 montre que $\mathcal{F}$ est une fibration elliptique. Ceci conclut la preuve du théorème. □

## 4. Applications birationnelles préservant une fibration rationnelle

Le but de cette section est de démontrer la proposition suivante.

**Proposition 4.1.** *Soit $\varphi$ une transformation birationnelle d'une surface complexe compacte dont la croissance des degrés est linéaire. Alors $\varphi$ préserve une fibration rationnelle et tout feuilletage $\varphi$-invariant coïncide avec cette fibration.*

L'existence d'une fibration rationnelle invariante résulte de l'hypothèse sur la croissance des degrés (voir 2.2). Ce qu'il faut montrer, c'est qu'il n'y a pas d'autre feuilletage invariant. Avant de donner la preuve, nous commençons par quelques remarques concernant les germes de transformations birationnelles holomorphes et les germes de feuilletages invariant par de telles transformations.

4.1. **Étude locale des feuilletages invariants.** Soit $\varphi \in \mathrm{Bir}(\mathbb{C}^2, 0)$ un germe holomorphe induit par une transformation birationnelle au voisinage de l'origine dans $\mathbb{C}^2$. En tout point où $\varphi$ est un biholomorphisme local, sa différentielle est inversible. L'ensemble critique de $\varphi$ est donc égal au diviseur contracté par $\varphi$. Nous supposons dans ce paragraphe que ce diviseur passe par 0 ; nous dirons alors que l'origine est un *point critique fixe.*

D'après le théorème de factorisation des transformations birationnelles, il existe un voisinage $\mathcal{U}$ de l'origine, une contraction $\pi : \tilde{\mathcal{U}} \to \mathcal{U}$ et un difféomorphisme $\sigma : \mathcal{U} \to \hat{\mathcal{U}}$ ($\hat{\mathcal{U}} \subset \tilde{\mathcal{U}}$) tels que $\varphi = \pi \circ \sigma$ (voir la figure 1). Nous noterons $E$ le diviseur exceptionnel $\pi^{-1}\{0\}$. Puisque 0 est un point critique, $\sigma(0)$ appartient à $E$.

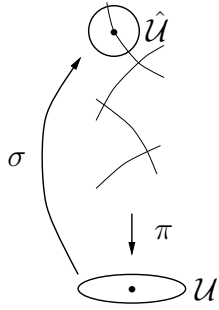

Figure 1. Factorisation locale.

Le diviseur $E$ étant obtenu par une suite d'éclatements successifs, c'est un arbre de courbes rationnelles ; en particulier, $E$ est à croisements normaux. L'ensemble analytique contracté par $\varphi$ coïncide avec $\sigma^{-1}(E)$ et possède donc une ou deux composantes irréductibles suivant que $\sigma(0)$ est un point lisse ou une singularité de $E$. Quitte à changer $\varphi$ en $\varphi \circ \varphi$, nous pouvons donc supposer que l'ensemble critique de $\varphi$ coïncide avec celui de ses itérés. L'ensemble critique est alors invariant par $\varphi$ et par $\varphi^{-1}$.

**Remarque 4.1.** Cette propriété de l'ensemble critique, essentielle dans la suite, a déjà été appliquée avec succès à l'étude des formes normales des points fixes attractifs d'applications birationnelles (voir [11]).

Soit $\mathcal{F}$ un germe de feuilletage holomorphe au voisinage de 0, éventuellement singulier à l'origine. Si $\mathcal{F}$ est $\varphi$-invariant, les feuilletages $\pi^*(\mathcal{F})$ et $\sigma_*(\mathcal{F})$ coïncident sur $\hat{\mathcal{U}}$. Après une suite d'éclatements, il existe donc un point du diviseur exceptionnel au voisinage duquel le feuilletage éclaté est difféomorphe au feuilletage initial. Nous allons voir que cette propriété est très restrictive.

**Proposition 4.2.** *Soit $\varphi$ un germe de transformation birationnelle au voisinage de l'origine dans $\mathbb{C}^2$ admettant l'origine comme point critique fixe. Si $\mathcal{F}$ est un germe de feuilletage holomorphe invariant par $\varphi$, alors :*

(i) *l'ensemble critique de $\varphi$ est tangent à $\mathcal{F}$ ;*

(ii) *le feuilletage $\mathcal{F}$ est lisse ou présente une singularité réduite non dégénérée à l'origine. Dans ce dernier cas, l'union des deux séparatrices de $\mathcal{F}$ coïncide avec l'ensemble critique de $\varphi$.*

*Démonstration.* Le théorème de désingularisation de Seidenberg affirme qu'après un nombre fini d'éclatements les singularités du feuilletage $\mathcal{F}$ sont réduites et qu'elles le restent si nous continuons à éclater. Dans notre situation, une singularité isomorphe à la singularité initiale de $\mathcal{F}$ apparaît après un nombre non nul d'éclatements : le feuilletage $\mathcal{F}$ est donc nécessairement réduit. En particulier, il ne passe qu'un nombre fini de courbes intégrales de $\mathcal{F}$ par l'origine. La courbe contractée par $\varphi$ ne peut donc pas être transverse au feuilletage, ce qui montre le point (i).

Supposons maintenant que $\mathcal{F}$ possède une singularité à l'origine qui soit de type selle-noeud. Le feuilletage $\mathcal{F}$ est formellement équivalent au feuilletage donné par la forme

$$\omega = x(1 + \lambda y^p)dy - y^{p+1}dx, \tag{4.1}$$

où $\lambda$ est un nombre complexe que nous appellerons "*invariant de Dulac*" et $p$ est un entier strictement positif (son nombre de Milnor). Ces deux nombres sont des invariants formels du feuilletage ; en particulier, ce sont des invariants analytiques. Lorsqu'on effectue un éclatement au niveau d'une singularité selle-noeud, il apparaît deux singularités sur le diviseur exceptionnel : elles se situent aux points d'intersection avec les transformées strictes de la séparatrice forte et de la séparatrice faible du feuilletage $\mathcal{F}$. Au niveau de la séparatrice forte, la singularité est non dégénérée ; si l'on éclate à nouveau ce point, il n'apparaît donc jamais de singularité selle-noeud. Au niveau de la séparatrice faible, la singularité est de type selle-noeud et est donnée formellement par

$$\tilde{\omega} = x(1 + (\lambda - 1)y^p)dy - y^{p+1}dx. \tag{4.2}$$

En ce point, la séparatrice faible coïncide bien sûr avec la transformée stricte de celle de $\mathcal{F}$. En résumé, après un nombre fini d'éclatements, une seule singularité est de type selle-noeud, elle se situe au niveau de la transformée stricte de la séparatrice faible et son invariant de Dulac a chuté d'au moins une unité. On ne peut donc jamais obtenir une singularité difféomorphe à la singularité initiale. En revenant aux hypothèses de la proposition, ceci montre que le feuilletage $\mathcal{F}$ présente un point régulier ou une singularité réduite non dégénérée à l'origine.

Supposons dorénavant que $\mathcal{F}$ possède une singularité réduite non dégénérée ; le feuilletage est donc donné au voisinage de 0 par une forme du type

$$\omega = x\, dy + \beta y\, dx + O(2), \tag{4.3}$$

où $\beta$ est un nombre complexe non nul. Il reste à montrer que les deux séparatrices de $\mathcal{F}$ à l'origine sont contenues dans l'ensemble critique de $\varphi$ ou, ce qui revient au même, que la singularité $\sigma(0)$ de $\pi^*(\mathcal{F})$ est un point singulier de $E$. Si ce n'est pas le cas, l'une des séparatrices de $\mathcal{F}$ n'est pas contenue dans l'ensemble critique ; nous la noterons $C$ et choisirons les coordonnées locales pour que $C$ soit tangent à l'axe $\{x = 0\}$. Le point $\sigma(0)$ doit être situé à l'intersection de la transformée stricte de $C$ avec $E$. En éclatant toujours la transformée stricte de $C$, nous avons donc obtenu une singularité isomorphe à la singularité initiale. Nous pouvons alors trouver des coordonnées $(u, v)$ au voisinage de $\sigma(0)$, telles que la courbe $\{u = 0\}$ soit incluse dans la transformée stricte de $C$ et que $\pi(u, v) = (uv^n, v)$. Au voisinage de $\sigma(0)$, le feuilletage $\pi^*(\mathcal{F})$ est donc défini par la forme

$$\hat{\omega} := \frac{\pi^*\omega}{v^n} = u \, dv + (\beta + n)v \, du + O(2). \tag{4.4}$$

Le difféomorphisme $\sigma$ envoie l'ensemble critique $\{y = 0\}$ sur le diviseur exceptionnel $E$, *i.e.* sur l'axe $\{v = 0\}$, et établit un isomorphisme entre les feuilletages donnés par ces deux formes différentielles. Les indices de Camacho-Sad satisfont donc $\beta = \beta + n$. L'entier $n$ devrait être nul, une contradiction. $\square$

4.2. **Preuve de la proposition 4.1.** Soit $\varphi : X \circlearrowleft$ une transformation birationnelle d'une surface complexe compacte préservant une fibration rationnelle $\pi : X \to B$. L'hypothèse sur la croissance des degrés de $\varphi$, à savoir $\|\varphi^{n*}\| \sim C^{ste} \cdot n$, est équivalente au fait que $\varphi$ n'est pas birationnellement conjuguée à un automorphisme (voir section 2). L'ensemble critique de $\varphi$, noté $\mathcal{C}(\varphi)$, et son ensemble d'indétermination $I(\varphi)$ sont donc tous deux non vides.

Supposons que $\varphi$ préserve un feuilletage holomorphe $\mathcal{F}$ distinct de la fibration $\pi$. Quitte à éclater suffisamment $X$, nous pouvons supposer que $\mathcal{F}$ est un feuilletage réduit. Il est alors possible d'effectuer de nouveaux éclatements pour rendre $\varphi$ algébriquement stable (voir section 2). Ceci signifie que, pour tout entier positif $n$, l'ensemble $\varphi^n(\mathcal{C}(\varphi) \setminus I(\varphi^n))$ est disjoint de $I(\varphi)$. Nous supposerons donc dans toute la suite que $\mathcal{F}$ est réduit et que $\varphi$ est algébriquement stable. Par contre, la nouvelle fibration rationnelle, toujours notée $\pi$, peut posséder des fibres singulières. Par construction, ces fibres sont des arbres de courbes rationnelles d'auto-intersection strictement négative.

**Remarque 4.2.** Puisque $\varphi$ préserve la fibration $\pi : X \to B$ , $\mathcal{C}(\varphi)$ est une union de courbes rationnelles lisses incluses dans des fibres de $\pi$.

**Remarque 4.3.** Puisque $\mathcal{F}$ est réduit, il ne possède pas de singularité dicritique. Les courbes contractées par $\varphi$ sont donc $\mathcal{F}$-invariantes.

**Lemme 4.3.** *Il existe un morphisme birationnel $\epsilon : X \to X'$ de telle sorte que*

(i) *l'application birationnelle $\epsilon \circ \varphi \circ \epsilon^{-1}$ soit algébriquement stable et $\epsilon_*\mathcal{F}$ soit un feuilletage réduit,*

(ii) *les fibres singulières de la fibration $\pi \circ \epsilon^{-1}$ ne contiennent pas de composante de l'ensemble critique de $\epsilon \circ \varphi \circ \epsilon^{-1}$.*

*Démonstration.* Choisissons une fibre singulière $F$ de la fibration contenant une composante de l'ensemble critique de $\varphi$. L'intersection $\mathcal{C}(\varphi) \cap F = V$ est une union de courbes rationnelles d'auto-intersection strictement négative. Factorisons $\varphi$ en une suite d'éclatements de points $\pi_1$ suivie d'une suite de contractions de courbes exceptionnelles de première espèce $\pi_2$. Notons $V'$ la transformée stricte de $V$ par $\pi_1$.

Si $V'$ ne possède aucune courbe d'auto-intersection $-1$, le dernier éclatement de $\pi_1$ peut être aussi choisi comme la première contraction de $\pi_2$. Quitte à choisir $\pi_1$ et $\pi_2$ minimales pour leur nombre d'éclatements, nous pouvons donc supposer que $V'$ contient une courbe $E'$ d'auto-intersection $-1$. Notons $E = \pi_1(E')$. Puisque $E$ est une courbe d'auto-intersection strictement négative, son auto-intersection est elle-même égale à $-1$ et aucun centre d'éclatements de $\pi_1$ n'est sur $E$.

Nous pouvons donc contracter la courbe $E$ sur un point $p$ et étudier le feuilletage induit par $\mathcal{F}$ en $p$. Le morphisme $\pi_1$ induit un isomorphisme d'un voisinage de $E'$ sur $E$. Le morphisme $\pi_2$ se décompose en la contraction de $E'$ sur un point soit $\pi_3$ suivi d'une suite de contractions $\pi_4$ ; ainsi, $\pi_2 = \pi_4 \circ \pi_3$. Au voisinage de $p$, le feuilletage induit par $\mathcal{F}$ est donc isomorphe à $\pi_4^* \mathcal{F}$ au voisinage de $r := \pi_3(E')$. Comme $\mathcal{F}$ est réduit, $\pi_4^* \mathcal{F}$ aussi, c'est donc aussi le cas pour le feuilletage induit par $\mathcal{F}$ en $p$.

Enfin, l'application induite par $\varphi$ reste algébriquement stable, car $I(\varphi) \cap E$ est inclus dans les centres d'éclatements de $\pi_1$ situés sur $E$ et est donc vide par ce qui précède. L'application induite par $\varphi$ en $p$ est donc holomorphe.

Nous avons montré que la contraction de $E$ n'altère pas les propriétés de $\mathcal{F}$ et $\varphi$. Quitte à itérer ce procédé nous pouvons donc supposer que les fibres singulières n'intersectent ni l'ensemble critique de $\varphi$ ni celui des itérés $\varphi^n$. La composée des contractions nécessaires est le morphisme $\epsilon$ cherché. □

**Lemme 4.4.** *Quitte à prendre un itéré de $\varphi$, toute composante critique $V$ de $\varphi$ est contractée par $\varphi$ sur un point $p \in V$ qui est fixe et par $\varphi^{-1}$ sur un point $q \in V$ distinct de $p$.*

*Démonstration.* Soit $V$ une composante critique de $\varphi$. D'après la remarque 4.3, $V$ est nécessairement invariante par $\mathcal{F}$. Si $\mathcal{F}$ possède une infinité de courbes rationnelles invariantes, un théorème de J.P. Jouanolou (voir [18] par exemple) montre que le feuilletage est induit par une fibration rationnelle. Dans ce cas, la preuve de la proposition 4.1 est terminée. Sinon l'ensemble des fibres de la fibration rationnelle qui sont aussi invariantes par $\mathcal{F}$ est un ensemble fini $\varphi$-invariant. Quitte à prendre un itéré de $\varphi$ nous pouvons supposer que chacune de ces fibres est invariante. Comme $\varphi$ est algébriquement stable, une fibre critique $V$ est contractée sur un point $p \in V$ qui n'est pas d'indétermination. Le même argument s'applique pour $\varphi^{-1}$. □

*Démonstration de la proposition 4.1.* Appliquons les lemmes 4.3 et 4.4. Quitte à faire un changement de variable birationnel, toute composante critique $V$ est une courbe rationnelle lisse d'auto-intersection nulle qui, de plus, est $\varphi$-invariante et se contracte sur un point critique fixe $p \in V$. En appliquant la proposition 4.2, nous en déduisons que $V$ est $\mathcal{F}$-invariante et que le feuilletage est régulier en $p$. De même $\mathcal{F}$ est régulier au point $q = \varphi^{-1}(V) = I(\varphi) \cap V$. Si un point $r \in V \setminus \{p, q\}$ est une singularité de $\mathcal{F}$, celle-ci est de type non-dégénéré et possède une courbe intégrale transverse à $V$, ou de type selle-noeud et possède une courbe intégrale formelle transverse à $V$. L'image par $\varphi$ de cette courbe transverse est une courbe intégrale du feuilletage (éventuellement formelle), et par invariance elle reste transverse à $V$ et passe par $p$. Ceci est impossible. Le feuilletage $\mathcal{F}$ ne contient donc pas de singularité sur $V$. Le théorème de stabilité de Reeb montre alors que $\mathcal{F}$ coïncide avec la fibration. □

## 5. Feuilletages dont le groupe de symétries birationnelles est infini

Le but de cette section est de démontrer le théorème 1.2. Nous supposerons donc que le groupe $\text{Bir}(\mathcal{F})$ contient *strictement* $\text{Aut}(\mathcal{F})$ quelque soit le modèle birationnel du feuilletage $\mathcal{F}$ étudié. Suivant [6], nous montrerons que $\text{Bir}(\mathcal{F})$ est de type Riccati pour une fibration rationnelle $\pi$. Lorsque tous les éléments de $\text{Bir}(\mathcal{F})$ préservent cette fibration, ou bien $\mathcal{F}$ coïncide avec cette fibration, ou bien $\text{Bir}(\mathcal{F})$ est fini et est égal à $\text{Aut}(\mathcal{F})$ après un changement de variable birationnel ad-hoc (paragraphe 5.4). Sinon, nous montrerons que $\mathcal{F}$ est Riccati pour une deuxième fibration rationnelle et que la situation est conjuguée à celle des exemples 1.3.

### 5.1. Le cas des feuilletages linéaires de $\mathbb{P}^1 \times \mathbb{P}^1$. 
Avant de débuter la preuve du théorème, décrivons la situation lorsque le feuilletage $\mathcal{F}$ est le feuilletage linéaire de $\mathbb{P}^1 \times \mathbb{P}^1$ déterminé par la 1-forme holomorphe

$$\omega = w \, dz + \alpha z \, dw, \tag{5.1}$$

où $\alpha$ est un nombre complexe non nul.

Ce feuilletage induit une fibration rationnelle si, et seulement si $\alpha$ est un nombre rationnel. Dans ce cas, le groupe $\text{Bir}(\mathcal{F})$ contient strictement $\text{Aut}(\mathcal{F})$, et ceci reste vrai même si un changement de variable birationnel est effectué ; les transformations birationnelles algébriquement stables construites dans l'exemple 1.2 suffisent en effet à établir cette remarque.

Lorsque $\mathcal{F}$ n'induit pas une fibration rationnelle, les seules courbes algébriques invariantes par $\mathcal{F}$ sont les quatre courbes $\{zw = 0\} \cup \{zw = \infty\}$ de $\mathbb{P}^1 \times \mathbb{P}^1$ ; le complémentaire de ces courbes est isomorphe à $\mathbb{C}^* \times \mathbb{C}^*$. L'ensemble critique de tout élément $\varphi$ de $\mathrm{Bir}(\mathcal{F})$ est inclus dans le diviseur constitué des quatre courbes $\mathcal{F}$-invariantes car $\mathcal{F}$ n'a pas de singularité dicritique. Cette propriété étant aussi satisfaite par $\varphi^{-1}$, $\varphi$ induit un automorphisme monomial de $\mathbb{C}^* \times \mathbb{C}^*$ ; autrement dit, $\varphi(x,y) = (x^a y^b, x^c y^d)$ est associé à une matrice carrée dont les coefficients $(a,b,c,d)$ sont entiers et dont le déterminant $ad - bc$ vaut $\pm 1$. Le vecteur $(1, \alpha)$ (ou $(\alpha, 1)$) est alors un vecteur propre de cette matrice.

Pour que $\mathrm{Bir}(\mathcal{F})$ soit infini sans que $\mathcal{F}$ soit une fibration rationnelle, il faut donc qu'il existe un élément de $\mathrm{GL}(2, \mathbb{Z})$ pour lequel $(1, \alpha)$ est un vecteur propre et dont la trace a une valeur absolue strictement plus grande que 2. Dans ce cas, la transformation $\varphi$ qui lui est associée est une transformation birationnelle dont les degrés croissent exponentiellement ; cette transformation est algébriquement stable dès que les coefficients matriciels $a$, $b$, $c$ et $d$ sont positifs, et l'on peut toujours se ramener à ce cas après conjugaison dans $\mathrm{GL}(2, \mathbb{Z})$. Ainsi, le groupe $\mathrm{Bir}(\mathcal{F})$ contient toujours strictement $\mathrm{Aut}(\mathcal{F})$ même si le modèle birationnel de $\mathcal{F}$ est changé.

En conclusion, ou bien $\mathrm{Bir}(\mathcal{F})$ est fini, ou bien l'une des deux propriétés suivantes est satisfaite

(i) $\mathcal{F}$ est une fibration rationnelle ;
(ii) $(1, \alpha)$ est un vecteur propre d'une matrice inversible à coefficients entiers dont la trace est strictement plus grande que 2 en valeur absolue ;

dans ces deux cas $\mathrm{Bir}(\mathcal{F})$ est infini et contient strictement $\mathrm{Aut}(\mathcal{F})$ pour tout modèle birationnel de $\mathcal{F}$.

5.2. **Feuilletages de Riccati.** Un feuilletage $\mathcal{F}$ est un feuilletage de Riccati s'il existe une fibration rationnelle $\pi$ dont les fibres génériques n'ont pas de point de tangence avec $\mathcal{F}$. Une telle fibration est dite *adaptée* au feuilletage. Les singularités d'un feuilletage de Riccati sont de nature spéciale. Nous donnons sans preuve (*voir* par exemple [7]) le

**Lemme 5.1.** *Soit $\mathcal{F}$ un feuilletage de Riccati à singularités réduites, $\pi$ une fibration rationnelle adaptée à $\mathcal{F}$ et $V$ une fibre lisse de $\pi$ qui est $\mathcal{F}$-invariante. L'une des deux propriétés exclusives suivantes est alors vérifiée.*

(i) *$V$ possède deux singularités du feuilletage qui sont toutes deux non-dégénérées ou toutes deux de type selle-noeud avec $V$ comme séparatrice forte. Les indices de Camacho-Sad des séparatrices en ces deux points sont des nombres complexes opposés.*
(ii) *$V$ ne possède qu'une singularité, celle-ci est de type selle-noeud et $V$ est sa séparatrice faible.*

Le reste de cette section est consacré à la preuve du théorème 1.2.

## 5.3. $\mathcal{F}$ est un feuilletage de Riccati.
Commençons par faire quelques hypothèses sur le feuilletage $\mathcal{F}$ étudié. Quitte à effectuer un nombre fini d'éclatements, nous pouvons supposer que $\mathcal{F}$ est réduit et que les courbes rationnelles $\mathcal{F}$-invariantes sont lisses. Nous supposons de surcroît que le couple $(X, \mathcal{F})$ est minimal pour ces deux propriétés. L'union des courbes rationnelles $\mathcal{F}$-invariantes est noté $D_\mathcal{F}$ : puisque $\mathcal{F}$ est réduit, c'est un diviseur à croisements normaux. Nous supposerons enfin que $\mathcal{F}$ *ne définit pas une fibration rationnelle*.

Nous allons montrer que $\mathcal{F}$ est un feuilletage de Riccati. La preuve que nous donnons reprend les arguments de [6] tout en les précisant.

**Lemme 5.2.** *Si $\varphi$ appartient à Bir($\mathcal{F}$), l'ensemble critique $\mathcal{C}(\varphi)$ et le lieu d'indétermination $I(\varphi)$ sont inclus dans $D_\mathcal{F}$.*

*Démonstration.* Fixons un élément $\varphi$ de Bir($\mathcal{F}$) qui ne soit pas un automorphisme. Comme $\mathcal{F}$ est réduit, l'ensemble critique $\mathcal{C}(\varphi)$ est inclus dans $D_\mathcal{F}$ (voir la remarque 4.3). Soit $p$ un point d'indétermination de $\varphi$ situé en dehors de $D_\mathcal{F}$ (s'il en existe). Alors $p$ n'appartient pas à $\mathcal{C}(\varphi)$ et le diviseur $\varphi(p)$ ne contient donc aucun point d'indétermination de $\varphi^{-1}$ ; autrement dit, $\varphi^{-1}$ est un *morphisme birationnel* d'un voisinage du diviseur $\varphi(p)$ sur un voisinage de $p$. Le diviseur $\varphi(p)$ peut donc être contracté sur un point $q$ en conservant une surface lisse : l'application induite par $\varphi$ détermine un biholomorphisme local de $p$ sur $q$. Après contraction, le feuilletage $\mathcal{F}$ a donc la même structure au voisinage de $q$ qu'au voisinage de $p$ ; en particulier, il est réduit. Les composantes de $D_\mathcal{F}$ restent lisses car $p$ n'appartient pas à $D_\mathcal{F}$ (donc $\varphi(p)$ ne coupe pas $D_\mathcal{F}$). Ceci contredit la minimalité de $(X, \mathcal{F})$ et permet d'affirmer que l'ensemble d'indétermination $I(\varphi)$ est contenu dans le diviseur $D_\mathcal{F}$. □

Le but est maintenant de montrer le lemme 5.3 énoncé ci-dessous en poursuivant l'étude d'un élément $\varphi$ de Bir($\mathcal{F}$) qui n'est pas un automorphisme. Pour cela, nous utiliserons une factorisation minimale de $\varphi$ en la composée d'une suite d'éclatements $\pi_1$ et d'une suite de contractions $\pi_2$. La minimalité signifie que toute courbe de première espèce contractée par $\pi_1$ ne l'est pas par $\pi_2$ et vice-versa.

Choisissons un point $p$ de $I(\varphi)$ et une courbe exceptionnelle de première espèce $V$ contractée par $\pi_1$ sur $p$. Si celle-ci n'intersecte pas l'ensemble critique de $\pi_2$, son image par $\pi_2$ reste d'auto-intersection $-1$ et peut être contractée sur un point lisse. À nouveau, $\mathcal{F}$ resterait réduit et les composantes de $D_\mathcal{F}$ lisses. Donc $V$ intersecte $\mathcal{C}(\pi_2)$ et $\pi_2(V)$ est d'auto-intersection positive ou nulle. Nous avons donc montré que, lorsqu'on effectue successivement les contractions composant $\pi_2$, les images de $V$ ont une auto-intersection qui croît jusqu'à devenir positive ou nulle. Nous noterons $V_1$ la première image de $V$ d'auto-intersection nulle obtenue par ce procédé. Puisque $V_1$ est une courbe rationnelle lisse d'auto-intersection nulle, $V_1$ est une fibre d'une fibration rationnelle $\pi$.

**Remarque 5.1.** Les arguments qui viennent d'être présentés montrent que, si $p$ appartient à $I(\varphi)$, alors $\varphi^{-1}(p)$ contient une courbe d'auto-intersection positive.

Étudions la position du feuilletage par rapport à la fibration $\pi$ construite ci-dessus.

Comme le feuilletage $\mathcal{F}$ est réduit en $p = \pi_1(V)$, la courbe $V$ possède au plus deux singularités réduites. De plus, si l'une de ces singularités est de type selle-noeud $V$ en est la séparatrice forte et le feuilletage $\mathcal{F}$ présente également une singularité selle-noeud en $p$. Projetons la situation au voisinage de $V_1$ et notons $\mathcal{F}_1$ le feuilletage image de $\mathcal{F}$. Le feuilletage $\mathcal{F}_1$ possède au plus deux singularités sur $V_1$ et, si l'une d'entre elles est de type selle-noeud, alors $V_1$ en est la séparatrice forte en ce point. L'intersection entre $T_{\mathcal{F}_1}$ et $V_1$ se calcule à l'aide de la formule (2.3) :

$$[T_{\mathcal{F}_1}] \cdot [V_1] = 2 - \sum_{p \in V_1} Z(\mathcal{F}_1, V_1, p), \qquad (5.2)$$

où l'indice $Z(\mathcal{F}_1, V_1, p)$ est égal à 1 lorsque la singularité est non-dégénérée et lorsque $V_1$ est la séparatrice forte d'une singularité selle-noeud. Il s'ensuit que $[T_{\mathcal{F}_1}] \cdot [V_1] \in \{0, 1, 2\}$.

Choisissons maintenant une fibre générique $W$ de la fibration rationnelle $\pi$ induite par $V_1$. Puisque $W$ et $V_1$ sont homologues, $[T_{\mathcal{F}_1}] \cdot [W]$ est égal à $[T_{\mathcal{F}_1}] \cdot [V_1]$, en particulier, c'est un nombre positif. Si $\mathcal{F}_1$ ne coïncide pas avec la fibration, nous pouvons supposer que $W$ n'est pas $\mathcal{F}_1$-invariante. La formule (2.4) s'écrit

$$[T_{\mathcal{F}_1}] \cdot [W] = [W]^2 - \text{Tang}(\mathcal{F}_1, W) = -\text{Tang}(\mathcal{F}_1, W) \leq 0, \qquad (5.3)$$

ce qui entraîne $[T_{\mathcal{F}_1}] \cdot [W] = 0$, puis $\text{Tang}(\mathcal{F}, W) = 0$. Nous avons donc montré que $\mathcal{F}_1$ était Riccati pour $\pi$.

**Remarque 5.2.** Notons que toutes les fibres de $\pi$ ne sont pas nécesairement lisses et que nous avons perdu la propriété de minimalité imposée au début. Plus précisément, il se peut qu'il existe une courbe de première espèce se contractant sur un point réduit de $\mathcal{F}$. Nous pouvons toutefois supposer qu'une telle courbe n'est pas inclue dans une fibre de $\pi$ car, sinon, en la contractant la fibration $\pi$ est préservée et le feuilletage $\mathcal{F}$ reste réduit.

Dans la suite nous considèrerons exclusivement le feuilletage $\mathcal{F}_1$. Pour alléger les notations, ce feuilletage sera encore noté $\mathcal{F}$ et $V_1$ sera notée $V$. Nous disposons donc du

**Lemme 5.3.** *Quitte à faire un changement de variable birationnel, le feuilletage $\mathcal{F}$ est réduit de type Riccati pour une fibration rationnelle $\pi$. Celle-ci possède* une fibre $V$ qui est lisse, donc d'auto-intersection nulle, qui est $\mathcal{F}$-invariante et qui possède deux singularités du feuilletage. Cette courbe $V$ est contenue dans l'ensemble critique d'un élément de $\mathrm{Bir}(\mathcal{F})$. De plus, toute courbe de première espèce dont la contraction préserve le caractère réduit de $\mathcal{F}$ est transverse à la fibration.

**Remarque 5.3.** Comme $\mathcal{F}$ possède deux singularités le long de $V$, le diviseur $D_{\mathcal{F}}$ est une union de courbes incluses dans les fibres de $\pi$ et d'au plus deux courbes rationnelles transverses.

5.4. **Lorsque $\pi$ est invariante.** Dans ce paragraphe nous supposons en outre que, pour tout élément $\varphi$ de $\mathrm{Bir}(\mathcal{F})$, l'ensemble critique de $\varphi$ est une union de fibres de $\pi$. Il revient au même de supposer que tout élément de $\mathrm{Bir}(\mathcal{F})$ préserve $\pi$.

Si $\mathrm{Bir}(\mathcal{F})$ possède un élément d'ordre infini non birationnellement conjugué à un automorphisme, nous pouvons appliquer la proposition 4.1, ce qui contredit l'hypothèse suivant laquelle $\mathcal{F}$ n'est pas une fibration. Nous supposerons donc que tout élément de $\mathrm{Bir}(\mathcal{F})$ est birationnellement conjugué à un automorphisme. Nous allons montrer qu'il existe un morphisme birationnel qui permet de conjuguer globalement tous les éléments de $\mathrm{Bir}(\mathcal{F})$ à des automorphismes. Ceci contredira l'hypothèse de départ suivant laquelle $\mathrm{Aut}(\mathcal{F})$ est un sous-groupe stricte de $\mathrm{Bir}(\mathcal{F})$ pour tout modèle birationnel de $\mathcal{F}$.

Fixons un élément $\varphi$ de $\mathrm{Bir}(\mathcal{F})$ et une composante irréductible $C$ de $\mathcal{C}(\varphi)$. La courbe $C$ est contenue dans une fibre de $\pi$ et est contractée par $\varphi$ sur un point $p$. La remarque 5.1 montre que $\varphi^{-1}(p)$ contient une courbe d'auto-intersection positive ou nulle. Par suite $C$ est une fibre régulière de $\pi$.

**Remarque 5.4.** Ceci montre que toute courbe contractée par $\varphi$ a une auto-intersection nulle. En particulier, la section 4 de [12] montre que $\varphi$ est conjugué à un automorphisme après un nombre fini d'éclatements ; autrement dit, $\varphi$ est conjugué à un automorphisme par un *morphisme birationnel* $\epsilon : \hat{X} \to X$.

La fibre $C$ est contractée sur $p$ et $\varphi$ préserve la fibration. Il existe donc des points d'indétermination sur $C$ dont les images par $\varphi$ recouvrent la fibre passant par $p$. Notons $C_1$ cette fibre. Puisque $C_1$ contient une composante irréductible contractée par $\varphi^{-1}$, $C_1$ est une fibre lisse. Ainsi, $C$ contient un unique point d'indétermination de $\varphi$ (image de $C_1$ par $\varphi^{-1}$).

Si $\mathcal{F}$ ne possède qu'une singularité le long de $C$, celle-ci est de type selle-noeud et C en est la séparatrice faible (lemme 5.2). Comme $\phi$ contracte $C$, cela est impossible. Donc $\mathcal{F}$ possède deux singularités le long de C. En reproduisant le même raisonnement à partir de $C_1$, nous obtenons une suite de courbes $\mathcal{F}$-invariantes $C_0 = C$, $C_1 = \varphi(C_0)$, $C_2 = \varphi(C_1)$, ... Cette suite est nécessairement finie car elle est constituée de fibres le long desquelles $\mathcal{F}$ a deux singularités. Nous noterons $k+1$ le nombre de fibres ainsi construites, de sorte que $\varphi(C_k) = C_0$.

Soient $p_i$ et $q_i$ les singularités de $\mathcal{F}$ le long de la courbe $C_i$. Ces singularités sont réduites et possèdent donc des séparatrices (éventuellement formelles) transverses à la fibre $C_i$. Nous les noterons $\Delta_{p_i}$ et $\Delta_{q_i}$ et nous désignerons par $\nu_i$ l'indice de Camacho-Sad $\mathbf{CS}(\Delta_{p_i}, p_i)$. D'après le lemme 5.1,

$$\mathbf{CS}(\Delta_{q_i}, q_i) = -\nu_i. \tag{5.4}$$

D'après la remarque 5.4, une suite finie d'éclatements permet de conjuguer $\varphi$ à un automorphisme. Nous noterons $\varpi : Y \to X$ le morphisme déterminé par cette suite d'éclatements (a priori $Y$ et $\varpi$ dépendent de $\varphi$) et $\tilde{\varphi}$ l'automorphisme de $Y$ induit par $\varphi$. Soit $D_Y$ le diviseur $\varpi^{-1}(C_0 \cup ... \cup C_k)$. Nous noterons $\Delta'_{p_i}, \Delta'_{q_i}$ les transformées strictes de $\Delta_{p_i}$ et $\Delta_{q_i}$, et $p'_i, q'_i$ leur intersection respective avec $D_Y$. Quitte à échanger $p_i$ et $q_i$, nous pouvons supposer que l'automorphisme $\tilde{\varphi}$ envoie chaque séparatrice $\Delta'_{p_i}$ sur $\Delta'_{p_{i+1}}$ pour $i = 0, ..., k-1$, c'est-à-dire que $\tilde{\varphi}(p'_i) = p'_{i+1}$ pour $i = 0, ..., k-1$. De même $\tilde{\varphi}(q'_i) = q'_{i+1}$. Par contre, $\tilde{\varphi}(p'_k) = p'_0$ ou $q'_0$.

Lors d'un éclatement, les indices de Camacho-Sad chute d'une unité, donc $\mathbf{CS}(\Delta'_{p_i}, p'_i)$ et $\mathbf{CS}(\Delta_{p_i}, p_i)$ sont égaux modulo 1. Puisque l'automorphisme $\tilde{f}$ envoie $p_i$ sur $p_{i+1}$ nous obtenons

$$\nu_0 = \nu_1 = \ ... \ = \nu_k \text{ modulo } 1. \tag{5.5}$$

Réalisons une série de "flips" aux points $p_i$, $q_i$ en éclatant l'un de ces points, puis en contractant la transformée stricte de $C_i$ (voir figure 2). Cette opération remplace les indices de Camacho-Sad $\nu_i$ par $\nu_i + 1$ ou $\nu_i - 1$ suivant que l'on éclate $q_i$ ou $p_i$ respectivement. Ce procédé permet de se ramener au cas où $\nu_0 = \ ... \ = \nu_k$.

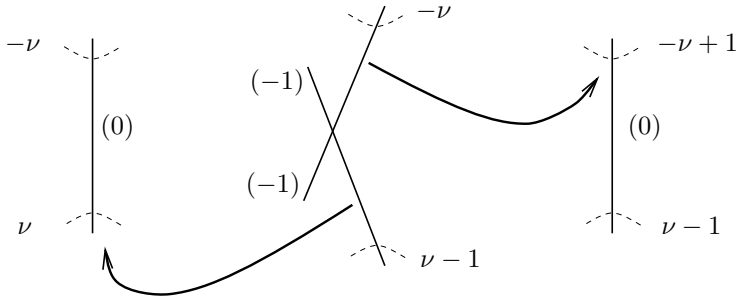

FIGURE 2. Un flip. Les indices de Camacho-Sad sont indiqués en lettres grecques et les auto-intersections sont entre parenthèses.

S'il existe une composante irréductible de $D_Y$ dont l'auto-intersection est égale à $-1$ et dont l'orbite sous $\tilde{\varphi}$ ne contient aucune transformée stricte $C'_i$ de $C_i$, toutes les courbes de cette orbite peuvent être contractées sans que $\tilde{\varphi}$ cesse d'être un automorphisme. Il s'ensuit qu'au moins un des $C'_i$ (soit $C'_0$) est d'auto-intersection égale à $-1$. En particulier, puisque $[C_0]^2 = 0$ les deux points $p_0$ et $q_0$ ne peuvent être simultanément des centres d'éclatements pour $\varpi$. Nous pouvons donc supposer que $p_0$ n'est pas un centre d'éclatement de $\varpi$. Le point $p'_0$ appartient donc à $C'_0$ et

$$\mathbf{CS}(\Delta'_{p_0}, p'_0) = \mathbf{CS}(\Delta_{p_0}, p_0). \tag{5.6}$$

L'image de $p'_0$ par l'automorphisme $\tilde{\varphi}$ est égale à $p'_1$, donc

$$\mathbf{CS}(\Delta'_{p_1}, p'_1) = \nu_0 = \nu_1 = \mathbf{CS}(\Delta_{p_1}, p_1). \tag{5.7}$$

Il s'ensuit que $p_1$ n'est pas un centre d'éclatement pour $\varpi$ et ainsi $\tilde{\varphi}(C'_0) = C'_1$. Par récurrence, on montre que $\tilde{\varphi}^i(C'_0) = C'_i$ pour $i = 1, ..., k$.

Supposons que $\tilde{\varphi}(p'_k) = p'_0$, ce qui assure $\tilde{\varphi}(C'_k) = C'_0$. Dans ce cas, l'union des composantes exceptionnelles de première espèce incluses dans $D_Y$ et distinctes de $C'_i$ peut être contractée sans créer de nouveau point d'indétermination pour $\tilde{\varphi}$, ce qui est absurde. Ainsi, $\tilde{\varphi}(p'_k) = q'_0$. Il s'ensuit que $-\nu_0 = \nu_0$ modulo 1. En réalisant alors une deuxième série de flips, nous pouvons maintenant imposer l'égalité

$$\nu_0 = \nu_1 = \ ... \ = \nu_k = -1/2 \text{ ou } 0. \tag{5.8}$$

Considérons tout d'abord le cas $\nu_i = -1/2$. Dans la surface $Y$, les transformées strictes $C'_i$ de $C_i$ sont toutes d'auto-intersection $-1$ et contiennent chacune le point $p'_i$. Par ailleurs, $\tilde{\varphi}(C'_i) = C'_{i+1}$ et $\tilde{\varphi}(p'_i) = p'_{i+1}$ pour $i = 0, \ ..., k-1$, ainsi que $\tilde{\varphi}(p'_k) = p'_0$. La courbe $E_0 = \tilde{\varphi}(C'_k)$ est donc d'auto-intersection $-1$ et contient $q'_0$. En ce point, l'indice de Camacho-Sad de $\Delta'_{q_0}$ est égal à $-1/2$. Lors de la suite d'éclatements déterminée par $\varpi$, la transformée stricte de $\Delta_{q_0}$ est éclatée exactement une fois, car sinon l'indice de Camacho-Sad serait strictement inférieur à $-1/2$. Par ailleurs, nous savons que la composante du diviseur exceptionnel $E_0$ contenant $q_0$ est d'auto-intersection $-1$ ; au-dessus de $q_0$ le morphisme $\varpi$ est donc exactement donné par l'éclatement simple de $q_0$. Nous avons donc prouvé que si l'on éclate une fois chacun des points $q_0, \ ..., q_k$, l'application $\varphi$ se relève en un automorphisme.

Dans le cas où $\nu_i = 0$, le même raisonnement montre que l'application $\varphi$ est un automorphisme sans avoir à éclater.

Pour conclure, remarquons maintenant que la suite de flips réalisée ne dépend pas du choix de $\varphi$ : elle dépend seulement de la nature des indices de Camacho-Sad aux points singuliers de $\mathcal{F}$. Nous avons donc montré l'existence d'un modèle birationnel de $\mathcal{F}$ pour lequel $\mathrm{Aut}(\mathcal{F})$ est égal à $\mathrm{Bir}(\mathcal{F})$. Ceci contredit nos hypothèses.

Nous pouvons désormais supposer que $D_\mathcal{F}$ contient une composante irréductible qui est transverse à la fibration $\pi$ et qui est contractée par un élément de $\mathrm{Bir}(\mathcal{F})$.

5.5. **Conjugaison à la situation monomiale.** Nous sommes dans la situation suivante : $\mathcal{F}$ est un feuilletage réduit de type Riccati (pour une fibration éventuellement singulière) qui possède une fibre lisse $V$ invariante par $\mathcal{F}$ et contractée par un élément de $\mathrm{Bir}(\mathcal{F})$ (lemme 5.3). Nous pouvons en outre supposer que $D_\mathcal{F}$ possède au moins une courbe $W$ qui est transverse à la fibration et qui est contractée par un élément $\varphi \in \mathrm{Bir}(\mathcal{F})$. Nous ne supposerons pas que cette application $\varphi$ est algébriquement stable.

Chaque point d'intersection de $V$ avec $W$ est un point singulier du feuilletage. Puisqu'une fibre lisse et invariante d'un feuilletage de Riccati possède au plus deux point singuliers, deux cas distincts apparaissent suivant que $V \cap W$ est constitué de un ou deux points.

**Lemme 5.4.** *La fibre $V$ possède deux singularités non-dégénérées d'indice non-rationnels. De même, la courbe $W$ possède exactement deux singularités non-dégénérées d'indice non-rationnels. L'une d'entre elle est située à l'intersection avec $V$.*

*Démonstration.* Soit $p$ un point de $V \cap W$. Si $p$ est une singularité de type selle-noeud alors $W$ en est la séparatrice forte car elle est contractée par un élément $\varphi$ de $\mathrm{Bir}(\mathcal{F})$. De même, $V$ est contractée par un élément $\varphi'$ de $\mathrm{Bir}(\mathcal{F})$ et doit donc être la séparatrice forte en ce point. Ceci est impossible. Ces singularités sont donc non-dégénérées. De plus, la somme des indices de Camacho-Sad des singularités sur $V$ est nul, ils ne peuvent donc être rationnels car $\mathcal{F}$ est réduit.

La formule de Camacho-Sad, appliquée à $W$ cette fois, montre que $W$ possède au moins deux singularités non-dégénérées d'indice non-rationnels. De telles singularités persistent après éclatement ou contraction car les indices de Camacho-Sad sont modifiés d'une unité par éclatement ou contraction. Ainsi, comme $W$ est contractée par $\varphi$, le point image $\varphi(W)$ est de type non-dégénéré d'indice non rationnel et toutes les composantes exceptionnelles de $\pi_2^{-1}(\varphi(W))$ possèdent au plus deux singularités non-dégénérées d'indice non-rationnels. C'est en particulier le cas pour la transformée stricte de $W$ par $\pi_1$, donc aussi pour $W$. Le lemme est démontré. □

Il se pourrait que $W$ ait aussi des singularités non dégénérées avec des indices rationnels. Excluons cette éventualité.

**Lemme 5.5.** *La courbe $W$ ne possède que deux singularités du feuilletage $\mathcal{F}$. Celles-ci sont toutes les deux non-dégénérées d'indice non-rationnels.*

*Démonstration.* Notons $p, q$ les deux singularités non-dégénérées de $W$ fournies par le lemme 5.4 et notons $\overline{W}$ la transformée stricte de $W$ par $\pi_1$. Décomposons $\pi_2$ en deux suites de contractions $\pi_2 = \widehat{\varpi} \circ \varpi$, de telle sorte que l'image de $\overline{W}$ par $\varpi$ soit d'auto-intersection $-1$. Nous noterons $W_1$ la courbe $\varpi(\overline{W})$ ; $W_1$ possède exactement deux singularités non-dégénérées d'indice non-rationnels que nous noterons $p_1, q_1$.

Regardons l'arbre des composantes exceptionnelles de $\varpi$ se contractant sur un point de $W_1$ distinct de $p_1$ et $q_1$. Si cet arbre est non-vide, il contient des composantes d'auto-intersection $-1$. Celles-ci ne peuvent être à la fois contractées par $\pi_1$ et $\pi_2$. Donc l'image par $\pi_1$ d'une telle courbe est une courbe ne possèdant par construction que des singularités de $\mathcal{F}$ d'indices rationnels non nuls. Une telle courbe sera notée $\widehat{C}$, et $C$ désignera son image par $\pi_1$.

Si $C$ n'est pas incluse dans une fibre de $\pi$, elle intersecte $V$ et possède donc une singularité non-dégénérée d'indice non-rationnel. Ceci est impossible et $C$ est donc incluse dans une fibre. Si son auto-intersection est négative, elle est égale à $-1$ et $\pi_1$ induit un isomorphisme de $\widehat{C}$ sur $C$. En particulier, nous pouvons contracter $C$ sur un point en conservant un feuilletage $\mathcal{F}$ réduit ce qui contredit les hypothèse faites.

Nous pouvons donc supposer que l'auto-intersection de $C$ est positive ou nulle ; il s'agit donc d'une fibre lisse de $\pi$ d'auto-intersection nulle. Elle intersecte $W$ et leur point d'intersection est une singularité réduite d'indice rationnel non nul. Le lemme 5.1 montre l'existence d'une singularité d'indice rationnel positif sur $C$. Ceci contredit le caractère réduit de $\mathcal{F}$.

Nous avons donc montré que $\overline{W}$ ne possède que deux singularités, non-dégénérées ; ces deux singularités sont projetées par $\pi_1$ sur les deux singularités d'indices non rationnels de $W$, et aucune autre singularité n'est créée par contraction. □

Dorénavant, nous distinguerons deux cas suivant que $V \cap W$ possède un ou deux points. Ces deux cas apparaissent respectivement dans les exemples 1.3-a et 1.3-b.

• *Les deux courbes $V$ et $W$ s'intersectent en un unique point $p$.*

En général, l'auto-intersection de $W$ n'est pas nulle, mais nous allons montrer que, quitte à faire des transformations birationnelles préservant la fibration, il est possible de se ramener à la situation $[V]^2 = [W]^2 = 0$.

Par le lemme 5.4, le feuilletage $\mathcal{F}$ possède deux singularités non-dégénérées d'indice non-rationnels sur $V$, soient $p := V \cap W$ et $q$. Fixons un petit disque $\Delta$ situé autour de l'image de $V$ par la projection $\pi$. Notons $(z, w)$ les coordonnées de $\pi^{-1}(\Delta)$ (identifié à $\Delta \times \mathbb{P}^1$) ; $\mathcal{F}$ est localement donné par la 1-forme holomorphe $\lambda w\, dz - z\, dw$ pour un nombre complexe irrationnel $\lambda$ convenable, les indices de Camacho-Sad en les deux singularités valent respectivement $\lambda$ et $-\lambda$, la fibration $\pi$ est donnée par $\pi(z, w) = z$ et $V$ coïncide avec la fibre $\{z = 0\}$ (voir par exemple [7] p.56).

Si $[W]^2$ est strictement positif, nous pouvons réduire l'auto-intersection de $W$ d'une unité en réalisant un flip en $p$ (en éclatant $p$ puis en contractant la transformée stricte de $V$). Si $[W]^2$ est strictement négatif, la même opération peut-être réalisée au point $q$ afin d'augmenter $[W]^2$ d'une unité. Dans les deux cas, le nouveau feuilletage est donné par la forme $(\lambda + \varepsilon)w\, dz - z\, dw$, où $\varepsilon$ est le signe de $[W]^2$, et possède encore deux singularités le long de la fibre $\{z = 0\}$. Celles-ci ont des indices non-rationnels et sont donc nécessairement réduites. Nous pouvons donc itérer ce procédé et, après une succession finie de flips, l'auto-intersection de $W$ nulle.

La courbe $W$ induit alors une fibration rationnelle $\pi'$. Le lemme 5.5 et la formule (2.4) montrent alors que $[T_\mathcal{F}] \cdot [W] = 2 - 2 = 0$. Le feuilletage $\mathcal{F}$ est donc également un feuilletage de Riccati vis-à-vis de la fibration $\pi'$. Puisque $[V] \cdot [W] = 1$, l'application $X \to \mathbb{P}^1 \times \mathbb{P}^1$ donnée par $\Psi(x) := (\pi(x), \pi'(x))$ est un morphisme birationnel. Notons $\mathcal{G}$ le feuilletage image de $\mathcal{F}$ par $\Psi$.

Comme le feuilletage $\mathcal{F}$ possède au plus deux singularités le long de $W$ (resp. $W$), la fibration $\pi$ (resp. $\pi'$) possède au plus une fibre totalement invariante distincte de $V$ (resp. $W$). Nous pouvons toujours supposer que $\Psi(V) = \{z = 0\}$, $\Psi(W) = \{w = 0\}$ et que les fibres totalement invariantes (si elles existent) sont envoyées sur $\{z = \infty\}$ et $\{w = \infty\}$. Il se pourrait qu'une fibre de $\pi'$ (resp. $\pi$), comporte une composante $\mathcal{F}$-invariante $C$ qui ne coupe pas $V$ (resp. $W$). Dans ce cas, $C$ est nécessairement contenue dans une fibre de $\pi$ (resp. $\pi'$). En particulier, $\Psi(C)$ est réduit à un point.

Dans le complémentaire de $\{z = \infty\} \cup \{w = \infty\}$ i.e. $\mathbb{C}^2$, le feuilletage $\mathcal{G}$ est donc transverse à toutes les courbes $z = c$ et $w = c$ dès que $c \neq 0$. De plus $\{zw = 0\}$ est $\mathcal{G}$-invariant. Le feuilletage $\mathcal{G}$ est défini par une 1-forme $\omega = p_1(z,w)\, dz + p_2(z,w)\, dw$ où $p_1, p_2$ sont des polynômes en $z, w$. Comme $\{z = 0\}$ est invariant, $z$ divise $p_2$. Par ailleurs pour $c$ non nul, la courbe $\{z = c\}$ est transverse à $\mathcal{G}$ donc $p_2(c, w)$ ne s'annule pas dans $\mathbb{C}$. Il s'ensuit que $\omega = w\, dz + \alpha z\, dw$ pour un certain nombre complexe non nul $\alpha$. La discussion de la section 5.1 permet alors de conclure.

• *Les deux courbes $V$ et $W$ s'intersectent en deux points disjoints.*

Par le lemme 5.5, les singularités de $\mathcal{F}$ sur $W$ sont incluses dans $V \cap W$. Nous noterons $p_0$ et $p_\infty$ ces deux singularités. Calculons $[W]^2$ en appliquant la formule de Camacho-Sad. D'après le lemme 5.4, les deux singularités sur $V$ sont non-dégénérées, donc

$$
\begin{align}
[W]^2 &= \mathbf{CS}(p_0, W) + \mathbf{CS}(p_\infty, W) \tag{5.9} \\
&= \mathbf{CS}(p_0, V)^{-1} + \mathbf{CS}(p_\infty, V)^{-1} \tag{5.10} \\
&= [V]^2 \cdot \mathbf{CS}(p_0, V)^{-1} \cdot \mathbf{CS}(p_\infty, V)^{-1} \tag{5.11} \\
&= 0 \, . \tag{5.12}
\end{align}
$$

Notons comme précédemment $\pi'$ la fibration rationnelle associée à $W$ pour laquelle $\mathcal{F}$ est aussi de type Riccati. Nous allons maintenant construire un revêtement double $Y$ de $X$ qui ramène la situation au cas précédent. Pour cela, considérons la sous-variété (singulière) $Y \subset X \times W$ définie par l'équation $\pi(p) = \pi(q)$ si $(p,q) \in X \times W$. Comme $[V] \cdot [W] = 2$, la projection naturelle $\varpi : Y \to X$ est un revêtement double ramifié. Par ailleurs,

$$\varpi^{-1}(V) = V \times \{p_0\} \cup V \times \{p_\infty\} \tag{5.13}$$

est l'union de deux courbes rationnelles lisses disjointes et d'auto-intersection nulle.

Pour comprendre $\varpi^{-1}(W)$ il nous faut étudier les singularités de $Y$. Pour tout point $m$ de $X$, notons $\sigma(m)$ le point d'intersection des fibres de $\pi$ et $\pi'$ passant par $m$ et distinct de $m$. L'involution $\sigma$ est définie hors des composantes communes aux deux fibrations[1] ; en particulier elle est définie au voisinage de $V \cup W$. Au niveau de $V \cap W$, $\sigma$ permute $p_0$ et $p_\infty$.

En un point $m$ de $W$ distinct de $p_0$ et $p_\infty$, le feuilletage $\mathcal{F}$ est lisse et $W$ est la feuille de $\mathcal{F}$ passant par ce point. La fibre $V_m$ de $\pi$ passant par $m$ intersecte $W$ en un ou deux points. Supposons qu'elle soit tangente à l'ordre 2 à $W$ en $m$. Le feuilletage $\mathcal{F}$ possèderait alors une feuille (en l'occurence $W$) tangente à l'ordre deux à $V_m$. Ceci est impossible car $\mathcal{F}$ est un feuilletage de Riccati pour la fibration $\pi$. Ainsi, toutes les fibres de $\pi$ intersectent transversalement $W$ et donc aussi toutes les fibres de $\pi'$ qui sont proches de $W$ [2]. Nous noterons $W_\varepsilon$, avec $\varepsilon \ll 1$, ces fibres proches de $W$. Considérons maintenant $\widehat{Y}$ la désingularisée minimale[3] de $Y$. L'application naturelle $\widehat{\varpi} : \widehat{Y} \to X$ induit un revêtement double au voisinage de $V \cup W$. La préimage du cycle $V \cup W$ est une union de quatre courbes rationnelles lisses d'auto-intersection nulles, ce sont les transformées strictes des variétés de $Y$

$$\varpi^{-1}(V) = V \times \{p_0\} \cup V \times \{p_\infty\} ; \tag{5.14}$$
$$\varpi^{-1}(W) = \{(p,p),\ p \in W\} \cup \{(p,\sigma(p)),\ p \in W\} . \tag{5.15}$$

Une composante de l'une de ces deux variétés intersecte toute composante de l'autre en un unique point ($(p_0, p_\infty)$, $(p_\infty, p_\infty)$, $(p_\infty, p_0)$ ou $(p_0, p_0)$). La surface $\widehat{Y}$ est donc isomorphe à $\mathbb{P}^1 \times \mathbb{P}^1$ éclaté. L'involution $\rho$ de $\widehat{Y}$ permutant les deux feuillets de $\widehat{\varpi}$ préserve le cycle $\widehat{\varpi}^{-1}(V \cup W)$ et permute les deux composantes de $\widehat{\varpi}^{-1}(V)$ et $\widehat{\varpi}^{-1}(W)$ respectivement. Si l'on envoie $\widehat{\varpi}^{-1}(V)$ sur $z = 0, \infty$ et $\widehat{\varpi}^{-1}(V)$ sur $w = 0, \infty$, l'involution $\rho$ s'écrit alors $\rho(z,w) = (1/z, 1/w)$ quitte à conjuguer par un automorphisme linéaire. On a donc montré que $X$ était la désingularisée du quotient de $\mathbb{P}^1 \times \mathbb{P}^1$ par $(z,w) \to (1/z, 1/w)$.

---

[1] on peut montrer que $\sigma$ est défini partout

[2] en particulier, $\pi$ possède des fibres singulières

[3] on peut montrer que $\widehat{Y}$ est la normalisée de $Y$

Le feuilletage $\mathcal{G}$ de $\mathbb{P}^1 \times \mathbb{P}^1$ qui est conjugué à $\mathcal{F}$ est, comme précédemment, défini par une 1-forme $z\, dw + \alpha w\, dz$. Puisque $\mathcal{F}$ n'est pqs une fibration rationnelle, $\mathcal{G}$ non plus. Tout élément $\varphi$ dans $\mathrm{Bir}(\mathcal{F})$ se le relève en un élément $\psi$ de $\mathrm{Bir}(\mathcal{G})$ tel que $\widehat{\varpi} \circ \psi = \varphi \circ \widehat{\varpi}$. Ainsi, $\mathrm{Bir}(\mathcal{G})$ est infini et la discussion du paragraphe 5.1 montre que $\mathrm{Bir}(\mathcal{G})$ contient une application monomiale dont la croissance des degrés est exponentielle. Celle-ci passe au quotient et $\mathrm{Bir}(\mathcal{F}) \setminus \mathrm{Aut}(\mathcal{F})$ contient un élément d'ordre infini de croissance des degrés exponentielle aussi.

Ceci termine la preuve du théorème 1.2 .

## 6. Autres approches

La preuve que nous avons donnée du théorème 1.2 et de son corollaire 1.3 s'appuie essentiellement sur les propriétés du feuilletage $\mathcal{F}$ et n'utilise pas la classification de McQuillan des feuilletages de dimension de Kodaira nulle. Deux autres approches conduisent également au corollaire 1.3.

Soit $\varphi$ une application birationnelle qui préserve un feuilletage $\mathcal{F}$, qui n'est pas birationnellement conjuguée à un automorphisme et dont la croissance des degrés est exponentielle.

6.1. **En utilisant les résultats de McQuillan.** Comme pour la démonstration du théorème 3.1, il est possible d'établir que la dimension de Kodaira $\mathrm{kod}(\mathcal{F})$ est nulle. Le point-clé consiste à démontrer que $\varphi$ préserve la partie positive du fibré canonique du feuilletage. Nour remercions vivement M. Brunella de nous avoir communiqué l'argument suivant. Factorisons $\varphi$ en la composition d'une suite d'éclatements et d'une suite de contractions, $\varphi = \pi_2 \circ \pi_1^{-1}$. Puisque $\mathcal{F}$ ne peut être une fibration rationnelle, le théorème de Miyaoka permet de supposer que $T_{\mathcal{F}}^*$ est numériquement effectif. Écrivons alors la décomposition de Zariski des fibrés canoniques de $\mathcal{F}$ et $\mathcal{G} = \pi_1^* \mathcal{F} = \pi_2^* \mathcal{F}$ :

$$[T_{\mathcal{F}}^*] = [P] + [N] \quad \text{et} \quad [T_{\mathcal{G}}^*] = [\overline{P}] + [\overline{N}], \tag{6.1}$$

où $[P]$ et $[\overline{P}]$ sont des classes nef. Puisque $T_{\mathcal{G}}^* = \pi_1^* T_{\mathcal{F}}^* + D$ où $D$ est un diviseur contractible, nous en déduisons l'égalité $\overline{P} = \pi_1^* P = \pi_2^* P$. Ceci garantit $\varphi^* P = P$ puis $P = 0$ car la croissance des degrés de $\varphi$ est exponentielle et $[P]$ est une classe nef. Nous avons donc établi l'égalité $\mathrm{kod}\,(\mathcal{F}) = 0$ (il s'agit ici de la dimension de Kodaira numérique ; dans le cas présent, elle coïncide avec la dimension de Kodaira définie dans l'introduction).

Nous pouvons donc relever la situation sur un revêtement fini $g : Y \to X$ pour lequel $g^* \mathcal{F}$ est engendré par un champ de vecteurs. L'application $\varphi$ se relève en une application birationnelle $\psi : Y \circlearrowleft$ dont la croissance des degrés doit être exponentielle et qui n'est pas birationnellement conjuguée à un automorphisme. Il s'ensuit que $Y$ est une surface rationnelle et il est alors aisé de conclure que la situation est conjuguée à celle de l'exemple 1.3.

6.2. **Approche dynamique.** Cette fois-ci nous privilégions l'étude de l'application $\varphi$ et montrons que son ensemble critique $\mathcal{C}(\varphi)$ est de même nature que celui d'une application monomiale (quatre courbes rationnelles lisses d'auto-intersection nulles dont le graphe dual est un carré), ou de son quotient par l'involution $(1/z, 1/w)$ (deux courbes rationnelles lisses d'auto-intersection nulles s'intersectant transversalement en deux points distincts).

Pour cela, il est possible de montrer que $\mathcal{C}(\varphi)$ est totalement invariant par $\varphi$, puis d'établir que $\mathcal{C}(\varphi)$ contient un unique cycle de courbes rationnelles et, enfin, qu'il est en fait réduit à un cycle. Une étude fine des singularités du feuilletage situées sur le cycle montre alors que $\mathcal{C}(\varphi)$ possède deux ou quatre courbes. Les arguments de la section 5.5 montrent ensuite que la situation est conjuguée à celle de l'exemple 1.3.

## 7. Quelques corollaires

7.1. **Dimension de Kodaira.** Le résultat principal de [28] apparaît comme corollaire des théorèmes principaux obtenus ici. La preuve de Pereira et Sanchez est une preuve directe qui utilise les premiers pas de la classification de McQuillan des feuilletages.

**Corollaire 7.1.** *([6], [28]) Supposons que $\mathcal{F}$ soit un feuilletage réduit dont la dimension de Kodaira est égale à 2. Les deux groupes $\mathrm{Aut}(\mathcal{F})$ et $\mathrm{Bir}(\mathcal{F})$ coïncident et sont d'ordre finis.*

Notons que la dimension de Kodaira du fibré cotangent d'un feuilletage générique de degré plus grand que 2 dans $\mathbb{P}^2$ est égale à 2.

*Démonstration.* Si $\mathrm{Aut}(\mathcal{F})$ et $\mathrm{Bir}(\mathcal{F})$ ne coïncident pas, le feuilletage $\mathcal{F}$ est de type Riccati (voir [6] ou le lemme 5.3). Dans ce cas, la dimension de Kodaira de $T^*_{\mathcal{F}}$ est au plus égale à 1 (voir [7]). Si $\mathrm{Aut}(\mathcal{F})$ est infini, le théorème 1.1 montre que $\mathcal{F}$ est donné par un champ de vecteurs, et alors $\mathrm{kod}(T^*_{\mathcal{F}}) = 0$, ou qu'il est donné par une fibration elliptique ou rationnelle, et alors $\mathrm{kod}(T^*_{\mathcal{F}}) \in \{-\infty, 0, 1\}$ ou, enfin, que $\mathcal{F}$ est un feuilletage de Riccati ou un feuilletage turbulent, et alors $\mathrm{kod}(T^*_{\mathcal{F}}) \in \{-\infty, 0, 1\}$. □

**Corollaire 7.2.** *Supposons que $\mathcal{F}$ soit défini sur une surface rationnelle $X$ et que $\mathrm{Bir}(\mathcal{F})$ soit infini. Alors, si $\mathcal{F}$ n'est pas une fibration rationnelle, sa dimension de Kodaira est nulle.*

*Démonstration.* Sous ces hypothèses, soit la situation est conjuguée à celle des exemples 1.1 ou 1.3 pour lesquels il est facile de voir que $\mathrm{kod}(\mathcal{F}) = 0$, soit $\mathcal{F}$ est invariant par un champ de vecteurs. Comme $X$ est rationnelle, nous pouvons nous ramener au cas (v) de la proposition 3.8. Le feuilletage est donc aussi induit par un champ de vecteurs holomorphe. Il s'ensuit que la dimension de Kodaira de $\mathcal{F}$ est nulle. □

**Remarque 7.1.** Soit $B$ une courbe de genre supérieur ou égal à deux et $E$ une courbe elliptique. Fixons $\alpha$ et $\beta$ des 1-formes holomorphes non-nulles respectivement sur $B$ et sur $E$ ; notons de même leur relevé sur la surface produit $B \times E$. La dimension de Kodaira du feuilletage engendré par la forme $\alpha + \beta$ sur $B \times E$ est égale à 1 et le feuilletage est invariant par tout champ de vecteurs parallèle à la fibration elliptique.

7.2. **Dynamique pseudo-Anosov.** Le corollaire suivant montre que les transformations rationnelles qui ont une entropie positive (il faut pour cela une croissance exponentielle des degrés) et qui préservent un feuilletage holomorphe sont linéarisables sur un ouvert de Zariski, préservent en fait deux feuilletages holomorphes et développent une dynamique de type Anosov ou pseudo-Anosov.

**Corollaire 7.3.** *Soit $\varphi$ une application birationnelle dont la croissance des degrés est exponentielle. Si $\varphi$ préserve un feuilletage, elle préserve deux feuilletages. Quitte à faire un revêtement ramifié, $\varphi$ est un automorphisme Anosov d'un tore complexe de dimension 2 ou une application monomiale de $\mathbb{P}^1 \times \mathbb{P}^1$.*

*Démonstration.* Il s'agit des exemples 1.1 ou 1.3. □

7.3. **Intégrabilité.** Concluons cette section d'applications en soulignant les liens étroits qui existent entre la classification obtenue et les propriétés d'intégrabilité des feuilletages. Nous remercions vivement J.V. Pereira pour ses explications concernant les propriétés d'intégrabilité des feuilletages. La démonstration du corollaire 7.4 proposée ci-dessous lui est entièrement dûe.

**Corollaire 7.4.** *Soit $\mathcal{F}$ un feuilletage d'une surface rationnelle possédant une infinité de symétries birationnelles. Alors $\mathcal{F}$ possède une intégrale première Liouvillienne.*

Une intégrale première rationnelle du feuilletage est une fonction rationnelle $h$ telle que la 1-forme $dh$ définisse le feuilletage. En d'autres termes, $\mathcal{F}$ admet une intégrale première rationnelle si et seulement si il définit une fibration. Nous renvoyons le lecteur à [30] pour les rudiments de la théorie des corps différentiels nécessaires à la définition précise des fonctions Liouvilliennes. Nous nous contenterons du critère d'intégrabilité de Singer ([30]) : un feuilletage algébrique $\mathcal{F}$ de $\mathbb{C}^2$ possède une intégrale première Liouvillienne si, et seulement s'il est défini par une 1-forme rationnelle $\omega$ telle que $d\omega = \eta \wedge \omega$ pour une 1-forme rationnelle $\eta$ *fermée*.

**Remarque 7.2.** Le critère de Singer signifie que le feuilletage est transversalement affine (sur un ouvert de Zariski). Dans [16], Ghys montre que les feuilletages stables et instables d'un difféomorphisme holomorphe Anosov (sur une surface complexe compacte) sont transversalement affines. Ceci suggère une autre approche pour obtenir une partie des résultats présentés dans cet article ; le lecteur notera cependant que nous n'avons fait aucune hypothèse sur la dynamique des symétries birationnelles (a priori, les symétries sont loin d'être Anosov).

*Démonstration du corollaire 7.4.* Supposons que le groupe Bir ($\mathcal{F}$) soit infini. Si le feuilletage $\mathcal{F}$ est une fibration rationnelle ou elliptique, il admet une intégrale première rationnelle. Sinon il existe un revêtement ramifié $\pi$ par $\mathbb{P}^1 \times \mathbb{P}^1$ ou par un tore $T$ tel que $\mathcal{F}$ se relève en un feuilletage défini par une 1-forme holomorphe fermée linéaire :

$$\omega_0 = z^{-1}\, dz + \lambda w^{-1}\, dw \ \text{ dans } \mathbb{P}^1 \times \mathbb{P}^1 \ \text{ ou } \ \omega_0 = dz + \lambda dw \ \text{ dans } T, \quad (7.1)$$

avec $\lambda \in \mathbb{C}^*$.

Choisissons une 1-forme rationnelle $\omega$ définissant $\mathcal{F}$. La forme $\pi^*\omega$ définit le même feuilletage que $\omega_0$ et nous pouvons donc trouver une fonction rationnelle $h$ telle que $h \cdot \pi^*\omega = \omega_0$. En différentiant cette égalité nous obtenons $d\pi^*\omega = h^{-1} \cdot dh \wedge \pi^*\omega$. Soit $\alpha := h^{-1} \cdot dh$. Le revêtement de $T$ ou $\mathbb{P}^1 \times \mathbb{P}^1$ sur $\mathbb{P}^1 \times \mathbb{P}^1$ induit une extension finie du corps des fonctions rationnelles de $\mathbb{P}^1 \times \mathbb{P}^1$. Notons que, dans notre cas, cette extension est galoisienne par construction. Nous pouvons définir $\pi_*\alpha$ comme la moyenne de toutes les images de $\alpha$ sous l'action du groupe de Galois de cette extension. En d'autres termes

$$\pi_*\alpha(x) := \frac{1}{d} \sum_{y \in \pi^{-1}(x)} (\pi_y^{-1})^*\alpha(y) \quad (7.2)$$

où $\pi_y^{-1}$ désigne la section de $\pi$ envoyant $x$ sur $y$ et $d$ est le degré du revêtement. Ceci définit une 1-forme rationnelle fermée $\eta$ de $\mathbb{P}^1 \times \mathbb{P}^1$ pour laquelle $d\omega = \eta \wedge \omega$. Ainsi, $\mathcal{F}$ admet une intégrale première Liouvillienne. □

## Appendix A. Automorphisme de surfaces rationnelles préservant un feuilletage.

Il s'agit de montrer la version précise suivante du corollaire 3.3.

**Corollaire A.1.** *Soit $T$ un tore complexe de dimension 2, muni d'un automorphisme $\varphi$ Anosov et d'un groupe fini d'automorphismes $G$ tel que $\varphi \circ G = G \circ \varphi$. Alors $T/G$ est un tore, une surface de Kummer, ou l'une des surfaces rationnelles suivantes :*

1. *La surface $T/G$ est isomorphe au quotient de $(\mathbb{C}/\mathbb{Z}[i])^2$ par le groupe d'ordre 4 engendré par l'homothétie de rapport $\sqrt{-1}$ et $\varphi$ est induit par l'action linéaire d'un élément de $SL(2, \mathbb{Z}[i])$ sur $\mathbb{C}^2$.*
2. *La surface $T/G$ est isomorphe au quotient de $(\mathbb{C}/\mathbb{Z}[j])^2$ par le groupe d'ordre 3 (resp. 6) engendré par l'homothétie de rapport $j$ (resp. $j^2$), et $\varphi$ est induit par l'action linéaire d'un élément de $SL(2, \mathbb{Z}[j])$.*

*Démonstration.* Nous allons décomposer la preuve en quatre étapes. Dans toute la suite, $\Gamma$ est un réseau de $\mathbb{C}^2$, $T$ est le tore $\mathbb{C}^2/\Gamma$ et $M$ est une transformation $\mathbb{C}$-linéaire de $\mathbb{C}^2$ qui est d'ordre fini $m$ et préserve $\Gamma$.

• **$M$ est d'ordre 2, 3, 4 ou 6.** En oubliant la structure complexe sur $\mathbb{C}^2$, nous pouvons supposer que $\Gamma$ est le réseau canonique $\mathbb{Z}^4$ et que $M$ correspond à un isomorphisme d'ordre fini $N \in GL(4, \mathbb{Z})$. Les valeurs propres de $N$ (donc de $M$) sont donc toutes des racines de l'unité dont le polynôme minimal sur $\mathbb{Z}$ est de degré au plus 4. Ainsi, en notant $\phi$ la fonction d'Euler,

$$\phi(m) \in \{1, 2, 3, 4\}. \tag{A.1}$$

Lorsque $\phi(m)$ vaut 1, $M$ est l'identité ; puisque $\phi$ ne prend jamais la valeur 3, nous obtenons $\phi(m) = 2$ ou 4.

Puisque $M$ est d'ordre fini, nous pouvons la diagonaliser ; nous noterons $\xi_1$ et $\xi_2$ les deux valeurs propres : il s'agit de deux racines de l'unité d'ordre $m$. Comme $\xi_i^{-1} = \overline{\xi_i}$, on a $M + M^{-1} = 2\,\mathrm{Re}\,(M)$. Si $\phi(m)$ est égal à 4, $m$ vaut 5, 8, 10, ou 12, et l'on constate sans peine que $\mathrm{Re}\,(\xi_1)$ et $\mathrm{Re}\,(\xi_2)$ sont des nombres réels irrationnels égaux ou opposés. Il s'ensuit que le groupe $H$ engendré par les matrices $\sum_0^{m-1} a_i M^i$ pour $a_i \in \mathbb{Z}$ n'est pas discret et adhère à la matrice nulle. Ceci contredit le fait que $H$ préserve le réseau $\Gamma$. Nous avons donc montré que $\phi(m)$ est égal à 2, ce qui assure que l'ordre de $M$ est égal à 2, 3, 4 ou 6.

• **$M$ est une homothétie.** Supposons maintenant que $M$ commute avec un automorphisme Anosov $\varphi$ de $T$. Si $M$ n'est pas une homothétie, elle admet deux vecteurs propres distincts $v_1$ et $v_2$, de valeurs propres $\xi_1$ et $\xi_2$, qui induisent deux plans réels $M$-invariants dans $\mathbb{R}^4$. Comme précédemment, $M$ induit une matrice $N$ préservant le réseau canonique $\mathbb{Z}^4$. Les deux plans invariants sont définis par $\ker[(N - \xi_i \,\mathrm{Id}\,)(N - \overline{\xi_i} \,\mathrm{Id}\,)]$. Puisque $M$ est d'ordre 2, 3, 4 ou 6, $\xi_i \in \{\pm 1, \pm i, \pm j, \pm j^2\}$, et $\mathbb{Q}[\xi_i] \cap \mathbb{R} = \mathbb{Q}$. Il s'ensuit que les plans invariants sont définis sur $\mathbb{Q}$ et que $\Gamma$ les intersecte.

Comme $\varphi$ commute à $M$, il préserve aussi ces plans. Mais $\varphi$ est Anosov et doit donc contracter strictement l'un de ces plans. Il ne peut donc préserver $\Gamma$, ce qui conduit à une contradiction. Nous avons donc montré que $M$ est une homothétie.

• **Le réseau $\Gamma$.** Nous pouvons désormais supposer que $M$ est une homothétie de rapport $\xi$, avec $\xi = i, j$ ou $-j$. Fixons un élément $\alpha$ dans le réseau $\Gamma$ dont la norme euclidienne est minimale. Au sein de la droite $\mathbb{C}\alpha$, les vecteurs $\{M^k \alpha, k \geq 0\}$ engendrent un sous-réseau du réseau $\Gamma \cap \mathbb{C}\alpha$. Par minimalité de $\alpha$, ces deux réseaux coïncident, en particulier, $\Gamma \cap \mathbb{C}\alpha$ est isomorphe à $\mathbb{Z}[\xi]$, c'est-à-dire à $\mathbb{Z}[i]$ ou $\mathbb{Z}[j]$.

Considérons maintenant le $\mathbb{C}$-espace vectoriel $\mathbb{C}^2/\mathbb{C}\alpha$ muni du réseau $\Gamma/(\Gamma \cap \mathbb{C}\alpha)$. Comme $M$ préserve la droite $\mathbb{C}\alpha$, elle induit une multiplication complexe par $\xi = i$ ou $j$ sur le réseau $\Gamma/(\Gamma \cap \mathbb{C}\alpha)$. Il s'ensuit que $\Gamma/(\Gamma \cap \mathbb{C}\alpha)$ est isomorphe à $\mathbb{Z}[i]$ ou $\mathbb{Z}[j]$ puis que $\Gamma$ lui-même est isomorphe à $\mathbb{Z}[i]^2$ ou $\mathbb{Z}[j]^2$.

• **Conclusion.** Nous revenons maintenant au problème initial. Nous disposons donc d'un groupe fini d'automorphismes du tore $T$, noté $G$, tel que $\varphi \circ G = G \circ \varphi$, où $\varphi$ est un automorphisme de type Anosov. Soit

$$\rho : G \to GL(2, \mathbb{C})$$

le morphisme envoyant un élément de $G$ sur sa partie linéaire. Notons $T(G)$ le noyau de $\rho$, i.e. le sous-groupe distingué constitué des éléments de $G$ de partie linéaire triviale. Alors $T(G)$ commute à $\varphi$ et nous pouvons passer au quotient par $T(G)$ ; nous obtenons un tore $T/T(G)$ sur lequel $\varphi$ induit à nouveau un automorphisme Anosov commutant au groupe fini $T/T(G)$. Ceci permet de supposer que $\rho$ est injectif. Si $M$ appartient à $G$, le commutateur $[M, \varphi]$ appartient à $G$ et $\rho[M, \varphi] = \mathrm{Id}$, donc $M$ et $\varphi$ commutent.

Prenons un élement $M$ de $G$ distinct de l'identité, donc avec une partie linéaire non triviale. Puisque $M$ commute à $\varphi$, l'automorphisme $M$ et l'automorphisme $\varphi$ possèdent un point fixe commun ; nous pouvons donc supposer que $M$ et $\varphi$ sont linéaires.

Les deux premières étapes montrent que $M$ est une homothétie $M = \xi \mathrm{Id}$ avec $\xi = -1, j, i, -j$ si $M$ est d'ordre 2, 3, 4 ou 6 respectivement. Le groupe $\rho(G)$ est donc isomorphe à l'un des quatre sous-groupes monogènes

$$\{\pm \mathrm{Id}\} \quad \cong \quad \mathbb{Z}/2, \tag{A.2}$$

$$\{j^k \mathrm{Id}\} \cong \mathbb{Z}/3 \quad \text{ou} \quad \{(-j)^k \mathrm{Id}\} \cong \mathbb{Z}/6, \tag{A.3}$$

$$\{i^k \mathrm{Id}\} \quad \cong \quad \mathbb{Z}/4. \tag{A.4}$$

Choisissons alors comme point base un point fixe commun à $\varphi$ et à un élement $M$ dont la partie linéaire engendre $\rho(G)$. Comme $\rho$ est injectif, tous les éléments de $G$ sont linéaires et $G$ est isomorphe à $\rho(G)$. Lorsque $G$ est d'ordre 2, la surface $T/G$ est une surface de Kummer. Lorque $G$ est d'ordre supérieur, la troisième étape montre que $\Gamma$ est le réseau $\mathbb{Z}[i]^2$ et $G$ le groupe $\langle i \mathrm{Id} \rangle$ ($G$ d'ordre 4), ou alors $\Gamma = \mathbb{Z}[j]^2$ et $G = \langle j \mathrm{Id} \rangle$ ($G$ d'ordre 3), ou alors $G = \langle -j \mathrm{Id} \rangle$ ($G$ d'ordre 6).

Ceci termine la preuve. □

## References

[1] Aldo Andreotti. Sopra le superficie algebriche che posseggono trasformazioni birazionali in sé. *Univ. Roma Ist. Naz. Alta Mat. Rend. Mat. e Appl.*, 9(5):255–279, 1950.

[2] Barth, Peters, and van de Ven. *Compact complex surfaces.* Springer-Verlag, 1984.

[3] Eric Bedford, Mikhail Lyubich, and John Smillie. Polynomial diffeomorphisms of $\mathbf{C}^2$. IV. The measure of maximal entropy and laminar currents. *Invent. Math.*, 112(1):77–125, 1993.

[4] Jean-Yves Briend and Julien Duval. Deux caractérisation de la mesure d'équilibre d'un endomorphisme de $\mathbf{P}^k(\mathbf{C})$. *Inst. Hautes Études Sci. Publ. Math.*, 93:145–159, 2001.

[5] Marco Brunella. Feuilletages holomorphes sur les surfaces complexes compactes. *Ann. Sci. École Norm. Sup. (4)*, 30(5):569–594, 1997.
[6] Marco Brunella. Minimal models of foliated algebraic surfaces. *Bull. Soc. Math. France*, 127(2):289–305, 1999.
[7] Marco Brunella. *Birational Geometry of Foliations.* Impa (voir http://www.impa.br/Publicacoes/Monografias/index.html), 2000.
[8] F. Campana and Th. Peternell. Cycle spaces. In *Several complex variables, VII*, pages 319–349. Springer, Berlin, 1994.
[9] Serge Cantat. Dynamique des automorphismes des surfaces projectives complexes. *C. R. Acad. Sci. Paris Sér. I Math.*, 328(10):901–906, 1999.
[10] Serge Cantat. Dynamique des automorphismes des surfaces K3. *Acta Mathematica*, 187(1):1–57, 2001.
[11] Charles Favre. Classification of 2-dimensional contracting rigid germs and Kato surfaces. I. *J. Math. Pures Appl. (9)*, 79(5):475–514, 2000.
[12] Charles Favre and Jeffrey Diller. Dynamics of bimeromorphic maps of surfaces. *Amer. J. Math.*, 123(6):1–29, 2001.
[13] Gerd Fischer. *Complex analytic geometry.* Springer-Verlag, Berlin, 1976. Lecture Notes in Mathematics, Vol. 538.
[14] John Erik Fornaess and Nessim Sibony. Complex dynamics in higher dimension. II. In *Modern methods in complex analysis (Princeton, NJ, 1992)*, pages 135–182. Princeton UNiv. Press, Princeton, NJ, 1995.
[15] T. Fujita. On Zariski problem. *Proc. Japan Acad. A*, 55:106–110, 1979.
[16] Étienne Ghys. Holomorphic Anosov systems. *Invent. Math.*, 119(3):585–614, 1995.
[17] Étienne Ghys. Feuilletages holomorphes de codimension un sur les espaces homogènes complexes. *Ann. Fac. Sci. Toulouse Math. (6)*, 5(3):493–519, 1996.
[18] Étienne Ghys. À propos d'un théorème de J.-P. Jouanolou concernant les feuilles fermées des feuilletages holomorphes. *Rend. Circ. Mat. Palermo (2)*, 49(1):175–180, 2000.
[19] Étienne Ghys and Alberto Verjovsky. Locally free holomorphic actions of the complex affine group. In *Geometric study of foliations (Tokyo, 1993)*, pages 201–217. World Sci. Publishing, River Edge, NJ, 1994.
[20] M. H. Gizatullin. Rational $G$-surfaces. *Izv. Akad. Nauk SSSR Ser. Mat.*, 44(1):110–144, 239, 1980.
[21] Mikhail Gromov. On the entropy of holomorphic map. *Manuscrit*, 1980.
[22] V. Guedj. Dynamics of polynomial maps of $\mathbf{C}^2$. preprint (2001), à paraître dans J. of the AMS.
[23] John H. Hubbard and Ralph W. Oberste-Vorth. Hénon mappings in the complex domain. I. *Inst. Hautes Etudes Sci. Publ. Math.*, 79:5–46, 1994.
[24] Henry B. Laufer. *Normal two-dimensional singularities.* Princeton University Press, 1971.
[25] David I. Lieberman. Compactness of the Chow scheme: applications to automorphisms and deformations of Kähler manifolds. In *Fonctions de plusieurs variables complexes, III (Sém. François Norguet, 1975–1977)*, pages 140–186. Springer, Berlin, 1978.
[26] Michael McQuillan. Diophantine approximations and foliations. *Inst. Hautes Études Sci. Publ. Math.*, 87:121–174, 1998.
[27] Luis Gustavo Mendes. Kodaira dimension of holomorphic singular foliations. *Bol. Soc. Brasil. Mat. (N.S.)*, 31:127–143, 2000.
[28] J.V. Pereira and P.F. Sanchez. Self bimeromorphisms of a foliation of general type. preprint.
[29] Nessim Sibony. Dynamique des applications rationnelles de $\mathbf{P}^k$. In *Dynamique complexe*, pages 97–185. Panoramas et Synthèses, S.M.F., 1999.
[30] M.F. Singer. Liouvillian first integral of differential equations. *Trans. of the AMS*, 1992.
[31] Yosi Yomdin. Volume growth and entropy. *Israel J. Math.*, 57(3):285–300, 1987.

Département de mathématiques, Université de Rennes, Rennes, France
*E-mail address*: `cantat@maths.univ-rennes1.fr`

Département de mathématiques, Université de Paris VII, Paris, France
*E-mail address*: `favre@math.jussieu.fr`